\documentclass[11pt, dina4,  amssymb]{article}
\usepackage{amssymb, amsmath,amsthm}
\usepackage{graphicx}
\usepackage{diagbox}
\usepackage{multirow}
\usepackage{epstopdf}
\usepackage{color}
\usepackage[margin=1in]{geometry}
\usepackage[normalem]{ulem}
\usepackage{xcolor}

\newcommand\redsout{\bgroup\markoverwith{\textcolor{red}{\rule[0.5ex]{2pt}{0.4pt}}}\ULon}

\newcommand{\bxi}{\textit{{\boldmath$\xi$}}}

\newcommand{\p}{\partial}
\newcommand{\og}{\omega}
\newcommand{\Og}{\Omega}
\newcommand{\fl}[2]{\frac{#1}{#2}}

\newcommand{\gm}{\gamma}
\newcommand{\nn}{\nonumber}
\newcommand{\ap}{\alpha}

\newcommand{\veps}{\varepsilon}

\newcommand{\Dt}{\Delta}

\newcommand{\be}{\begin{equation}}
\newcommand{\ee}{\end{equation}}
\newcommand{\ba}{\begin{array}}
\newcommand{\ea}{\end{array}}
\newcommand{\bea}{\begin{eqnarray}}
\newcommand{\eea}{\end{eqnarray}}
\newcommand{\beas}{\begin{eqnarray*}}
\newcommand{\eeas}{\end{eqnarray*}}
\newtheorem{remark}{Remark}[section]

\newcommand{\bx}{{\bf x} }
\newcommand{\bm}{{\bf m} }
\definecolor{ForestGreen}{rgb}{0.0, 0.5, 0.0}

\title{Highly accurate operator factorization methods for the integral fractional Laplacian and its generalization}
\author{
Yixuan Wu\thanks{Department of Mathematics and Statistics, Missouri University of Science and Technology, Rolla, MO 65409 (Email:  ywx7c@mst.edu)}, \ \
Yanzhi Zhang\thanks{Department of Mathematics and Statistics, Missouri University of Science and Technology, Rolla, MO 65409 (Email:  zhangyanz@mst.edu; URL:  {http://web.mst.edu/$\sim$zhangyanz})}}
\begin{document}
\date{}
\maketitle

\begin{abstract}
In this paper, we propose a new class of operator factorization methods to discretize the integral  fractional Laplacian $(-\Dt)^\fl{\ap}{2}$ for $\ap \in (0, 2)$. 
The main advantage of our method is to easily increase numerical accuracy by using high-degree Lagrange basis functions, but remain the scheme structure and computer implementation unchanged. 
Moreover, our discretization of the fractional Laplacian results in a symmetric (multilevel) Toeplitz differentiation matrix, which not only saves memory cost in simulations but  enables efficient computations via the fast Fourier transforms. 
The performance of our method in both approximating the fractional Laplacian and solving the fractional Poisson problems was detailedly examined. 
It shows that our method  has an optimal accuracy of ${\mathcal O}(h^2)$ for constant or linear basis functions, while ${\mathcal O}(h^4)$ if quadratic basis functions are used, with $h$ a small mesh size. 
Note that this accuracy holds for any $\ap \in (0, 2)$ and can be further increased if  higher-degree basis functions are used. 
If the solution of fractional Poisson problem satisfies $u \in C^{m, l}(\bar{\Og})$ for $m \in {\mathbb N}$ and $0 < l < 1$, then our method has an accuracy of ${\mathcal O}\big(h^{\min\{m+l,\, 2\}}\big)$  for constant and linear basis functions, while ${\mathcal O}\big(h^{\min\{m+l,\, 4\}}\big)$  for quadratic basis functions.
Additionally, our method can be readily applied to study  generalized fractional Laplacians with a symmetric kernel function, and numerical study on the tempered fractional Poisson problem demonstrates its efficiency. 
\end{abstract}


\section{Introduction}
\setcounter{equation}{0}
\label{section1}

The fractional Laplacian $(-\Dt)^\fl{\ap}{2}$, representing the infinitesimal generator of a symmetric $\ap$-stable L\'evy process, is a nonlocal generalization of the classical Laplace operator $-\Dt = -(\p_{xx} + \p_{yy} + \p_{zz})$. 
Over the recent decade,  it has been widely applied  to study anomalous diffusion  in many fields. 
However,  current understanding of the fractional Laplacian still remains limited in comparison to its classical counterpart. 
Specifically, the nonlocal nature of the fractional Laplacian introduces significant challenges in its theoretical analysis and numerical approximation. 
In this paper, we propose a new class of operator factorization methods to discretize the integral fractional Laplacian $(-\Dt)^\fl{\ap}{2}$ for $\ap \in (0, 2)$.  
One main advantage of our method is its flexibility to increase numerical accuracy by using high-degree basis functions. 

Let $\Og \subset {\mathbb R}^d$ (for $d = 1, 2, 3$) be an open bounded domain, and denote $\Og^c = {\mathbb R}^d\backslash\Og$ as the complement of $\Og$. 
The fractional Poisson problem with extended Dirichlet boundary condition takes the following form:
 \bea\label{Poisson}
(-\Dt)^\fl{\ap}{2}u(\bx) = f(\bx), &\ &\mbox{for} \ \ \bx \in\Og,\\
\label{BC}
u(\bx) = g(\bx), &&\mbox{for} \ \ \bx \in\Og^c. 
\eea
The fractional Laplacian $(-\Dt)^\fl{\ap}{2}$ is defined in a hypersingular integral form \cite{Landkof, Samko}:
\bea\label{fL-nD}
(-\Delta)^{\fl{\ap}{2}}u({\bf x}) 
= c_{d,\ap} \ {\rm P.V.}\int _{\mathbb{R}^d} \frac{u({\bf x})-u({\bf y})}{|{\bf x}-{\bf y}|^{d+\alpha}}\,d{\bf y}, \qquad \mbox{for \ $\alpha \in (0,2)$},
\eea
where P.V. stands for the Cauchy principal value,  $|\bx - {\bf y}|$ denotes the Euclidean distance between points $\bx$ and ${\bf y}$, and the normalization constant $c_{d,\ap}$ is defined as
\beas\label{DefC1ap}
 c_{d,\ap} = \fl{2^{\ap-1}\ap\,\Gamma\big(({d+\ap})/{2}\big)}{\sqrt{\pi^d}\,\Gamma\big(\displaystyle 1-{\ap}/{2}\big)}, \qquad  \mbox{for} \  \  \ d = 1, 2, 3
\eeas
with $\Gamma(\cdot)$ denoting the Gamma function.  
The operator \eqref{fL-nD} collapses to the identity operator as $\ap \to 0$, while it converges to the classical Laplacian $-\Dt$ as $\ap \to 2$. 
In the literature, the fractional Laplacian can be also defined via a pseudo-differential operator with symbol $|{\bf k}|^\alpha$  \cite{Landkof, Samko}, i.e., 
\begin{equation}
\label{pseudo}
(-\Delta)^{\fl{\alpha}{2}}u({\bx}) = \mathcal{F}^{-1}\big[|{\bf k}|^\alpha \mathcal{F}[u]\big], \qquad \mbox{for} \ \ \ap > 0,
\end{equation}
where $\mathcal{F}$ represents the Fourier transform, and $\mathcal{F}^{-1}$ denotes its inverse.  
It shows in \cite{Samko, Kwasnicki2017} that these two definitions are equivalent  on the Schwartz space. 
More discussion of the fractional Laplacian and its related nonlocal operators can be found in \cite{Du2012, Duo2019} and references therein.  
Note that the fractional Laplacian $(-\Dt)^\fl{\ap}{2}$ is rotational invariant, which is an important property in modeling isotropic anomalous diffusion especially when $d\geq 2$ \cite{Izsak2017}.

In the literature, numerical methods for the fractional Laplacian $(-\Dt)^\fl{\ap}{2}$ can be mainly classified in three groups based on which definition (e.g. \eqref{pseudo} or \eqref{fL-nD}) the method is developed on. 
The pseudo-differential definition in \eqref{pseudo} allows one to utilize the Fourier transform for both analysis and simulations. 
Based on this definition, it is natural to introduce the Fourier pseudospectral methods to discretize the fractional Laplacian $(-\Dt)^\fl{\ap}{2}$ on a bounded domain with periodic boundary conditions  \cite{Kirkpatrick2016}.  
This method can be applied for both classical ($\ap = 2$) and fractional ($\ap < 2$) Laplacians. 
However,  it is challenging to incorporate non-periodic boundary conditions into the pseudo-differential definition \eqref{pseudo}. 
For example, a meshfree pseudospectral method was proposed in \cite{Rosenfeld2019} to discretize $(-\Dt)^\fl{\ap}{2}$ with extended Dirichlet boundary conditions. 
In order to account Dirichlet boundary conditions, a large computational domain (much larger than physical domain $\Og$) as well as boundary approximations were introduced, and moreover numerical quadrature rules were required to approximate the Fourier integrals in definition (\ref{pseudo}).  
On the other hand, the integral definition in (\ref{fL-nD}) provides a pointwise formulation and thus can easily work with different boundary conditions.
Based on the integral definition, finite difference methods \cite{Huang2014, Duo2018, Duo-TFL2019} and finite element methods \cite{Acosta2017, Acosta2017-2D, Bonito2017,  Ainsworth2018B} have been developed to discretize the fractional Laplacian $(-\Dt)^\fl{\ap}{2}$ with extended Dirichlet boundary conditions. 
In contrast to the pseudo-differential definition, the integral definition in \eqref{fL-nD} is valid only for $\ap \in (0, 2)$, so do the resulting numerical methods \cite{Huang2014, Duo2018, Duo-TFL2019, Acosta2017, Acosta2017-2D, Bonito2017,  Ainsworth2018B}. 
Moreover, the above methods based on the integral definition are mainly limited to the second order accuracy.  
Recently, a new class of methods using both definitions \eqref{pseudo} and \eqref{fL-nD} of the fractional Laplacian $(-\Dt)^\fl{\ap}{2}$ were proposed  in  \cite{Burkardt0020, Wu0021}. 
These methods, based on radial basis functions, are meshfree and work for both classical and fractional Laplacians.

In this paper, we propose a novel  operator factorization method to discretize the integral  fractional Laplacian in \eqref{fL-nD}. 
Our method factorizes the integrand of $(-\Dt)^\fl{\ap}{2}$ in \eqref{fL-nD} as a product of central difference quotient function $\Phi_{d, \gm}(\bx, \bxi)$ and power function $\mu_{\gm}(\bxi)$, and then approximates  $\Phi_{d, \gm}(\bx, \bxi)$  in $\bxi$ using Lagrange basis functions $\varphi^p$ (for $p \in {\mathbb N}^0$). 
The main advantage of our method is that it can easily increase numerical accuracy by using high-degree Lagrange basis functions $\varphi^p$, but remain the scheme structure and computer implementation unchanged. 
Moreover,  our method results in a symmetric (multilevel) Toeplitz matrix, and thus algorithms via the fast Fourier transforms can be developed for their efficient simulations. 
Numerical studies show that our method with constant basis $\varphi^0$ or linear basis $\varphi^1$ has an optimal accuracy of ${\mathcal O}(h^2)$ with $h$ small mesh size, and this rate can be improved to ${\mathcal O}(h^4)$ if quadratic basis $\varphi^2$ is used. 
This accuracy can be further increased if high-degree basis functions are used. 
The performance of our method under different conditions are detailedly investigated and compared. 
If the solution of fractional Poisson problem satisfies $u \in C^{m, l}(\bar{\Og})$ for $m \in {\mathbb N}$ and $0 < l < 1$, then our method has an accuracy of ${\mathcal O}\big(h^{\min\{m+l,\, 2\}}\big)$  for constant and linear basis functions, while ${\mathcal O}\big(h^{\min\{m+l,\, 4\}}\big)$  for quadratic basis functions.
Our method can be readily applied to study the generalized fractional Laplacian with a symmetric kernel function, and numerical study on the tempered fractional Poisson problem is provided to demonstrate it. 

The paper is organized as follows. 
In Section \ref{section2}, we first introduce the operator factorization framework, and then derive our methods for both one- and two-dimensional cases.
Generalization to three dimensions is straightforward following the same lines as in two-dimensional scheme.  
In Section \ref{section3}, we examine the performance of our method in approximating the fractional Laplacian $(-\Dt)^\fl{\ap}{2}$ and  in solving fractional Poisson problems.  
Also, we generalize our method to solve the tempered fractional Poisson problems. 
Finally, conclusions are drawn in Section \ref{section4}. 

\section{Method of operator factorization}
\label{section2}
\setcounter{equation}{0}

Numerical methods for the  fractional Laplacian in (\ref{fL-nD}) still remain limited, and main challenges come from its nonlocality and strong singularity.  
So far, the existing finite difference/element methods are mostly the second-order accurate. 
In this section, we introduce a new class of methods based on operator factorization to discretize the  integral fractional Laplacian $(-\Dt)^\fl{\ap}{2}$ with extended Dirichlet boundary conditions. 
Our method has advantages of easily improving numerical accuracy without changing the scheme structure.

To introduce our method, we will first reformulate the integral fractional Laplacian (\ref{fL-nD}) into the form of operator factorization.  
For points $\bx$, ${\bf y} \in {\mathbb R}^d$, we define a vector $\bxi = \big(\xi^{(1)}, \, \xi^{(2)}, \, \cdots, \, \xi^{(d)}\big)$ with  $\xi^{(i)}=  |x^{(i)} - y^{(i)}|$ denoting the distance of $\bx$ and ${\bf y}$ in the $i$-th direction, and then rewrite the fractional Laplacian in (\ref{fL-nD}) as:
\bea \label{FL1}
(-\Dt)^{\fl{\ap}{2}}u(\bx) = -c_{d, \ap} \int_{{\mathbb R}_+^d} \bigg(\sum_{{\bf m} \in \varkappa_1} u(\bx + (-1)^\bm \circ \bxi) - 2^d u(\bx)\bigg)\, \fl{d\bxi}{|\bxi|^{d+\ap}},
\eea
where we denote ${\mathbb R}_+^d = [0, \infty)^d$, and ${\bf a} \circ {\bf b}$ represents the Hadamard product of ${\bf a}$ and ${\bf b}$.  
For $M \in {\mathbb N}$, the index set $\varkappa_M$ is defined as
\beas
\varkappa_M = \{(m_1, \, m_2, \, \cdots, \, m_d) \ | \ 0 \le m_i \le M, \ \mbox{for} \ 1 \le i \le d \}, 
\eeas
and the  vector  $(-1)^{\bf m} = \big((-1)^{m_1}, \, (-1)^{m_2}, \, \cdots, \, (-1)^{m_d} \big)$. 
Choosing a splitting parameter $\gm \in (\ap, 2]$, and introducing functions
\bea\label{eq0}
\Phi_{d, \gm}(\bx, \bxi) =\bigg(\sum_{{\bf m} \in \varkappa_1} u(\bx + (-1)^\bm \circ \bxi) - 2^d u(\bx)\bigg) {|\bxi|^{-\gm}}, \qquad \mu_\gm(\bxi) = |\bxi|^{\gm - (d + \ap)},
\eea
we then further rewrite  (\ref{FL1}) into the following operator factorization form:  
\bea\label{FL2}
(-\Dt)^{\fl{\ap}{2}}u(\bx) = -c_{d, \ap}\int_{{\mathbb R}_+^d}  \Phi_{d, \gm}(\bx, \bxi)\,\mu_\gm(\bxi) d\bxi.
\eea
In (\ref{FL2}), we factorize the integrand of the fractional Laplacian $(-\Dt)^\fl{\ap}{2}$ as a product of central difference quotient function $\Phi_{d, \gm}(\bx, \bxi)$ and power function $\mu_{\gm}(\bxi)$, which can be also viewed as a weighted integral of function $\Phi_{d, \gm}(\bx, \bxi)$. 
The above operator factorization was introduced in \cite{Duo2018, Duo-FDM2019} to develop finite difference methods for the integral fractional Laplacian $(-\Dt)^\fl{\ap}{2}$, and also  applied to solve the fractional Sch\"odinger equation in a box potential (with $\gm = 1+\ap/2$) \cite{Duo2015}. 
Note that even though both of them are based on the same operator factorization form, the method proposed in this work is significantly different from those in  \cite{Duo2018, Duo-FDM2019}. 
The  parameter $\gm$ plays an important role in determining the numerical accuracy of our method, and the optimal splitting parameter is $\gm = 2$ which leads to the smallest numerical errors and best accuracy rates. 
More discussion and illustrations will be given in Section \ref{section3}.

Denote the $d$-dimensional domain $\Og = (a_1,\, b_1) \times (a_2,\, b_2)  \cdots \times (a_d, b_d)$, and introduce 
\beas
\Upsilon = [0, \, L]^d \ \, \ \mbox{with}  \ \ L = \max_{1 \le i \le d}|b_i - a_i|; \qquad \Upsilon^c_+ = {\mathbb R}^d_+\backslash\Upsilon. 
\eeas
Then the integral in \eqref{FL2} can be divided  into two parts, i.e.
\bea \label{FL30}
(-\Dt)^{\fl{\ap}{2}}u(\bx) = -c_{d, \ap} \bigg(\int_{\Upsilon}  \Phi_{d, \gm}(\bx, \bxi) \mu_\gm(\bxi)\,d\bxi +
\int_{\Upsilon_+^c}  \Phi_{d, \gm}(\bx, \bxi) \mu_\gm(\bxi)\,d\bxi\bigg),  \qquad \mbox{for} \ \ \bx \in \Og.
\eea
Given the definition of $\Upsilon_+^c$, it is easy to conclude that for any $\bx \in \Og$ and  $\bxi \in \Upsilon_+^c$,\, the point $\big(\bx + (-1)^\bm \circ \bxi\big)$ for $\bm \in \varkappa_1$ locates outside of domain  $\Og$.
Noticing the extended Dirichlet boundary conditions in \eqref{BC}, we immediately obtain 
\bea\label{L2}
{\mathcal L}_{\Upsilon_+^c}^\ap u(\bx) &:=& \int_{\Upsilon_+^c}  \Phi_{d, \gm}(\bx, \bxi) \mu_\gm(\bxi)\,d\bxi\nn\\
&=& -2^du(\bx) \int_{\Upsilon^c_+} \fl{d\bxi}{|\bxi|^{d+\ap}}  + \sum_{{\bf m} \in \varkappa_1}\int_{\Upsilon^c_+} \fl{g\big(\bx + (-1)^\bm \circ \bxi\big)}{|\bxi|^{d+\ap}}\,d\bxi, \qquad \mbox{for} \ \ \bx \in \Og.\quad 
\eea
The above integrals are free of singularity, and thus can be accurately computed by traditional quadrature rules. 
Moreover, if homogeneous boundary conditions (i.e., $g(\bx) = 0$) are considered, the second term in \eqref{L2} vanishes.  

Next, we move to approximate the first integral in (\ref{FL30}), i.e.,
\bea\label{L1}
{\mathcal L}_{\Upsilon}^\ap u(\bx) &:=& \int_{\Upsilon}  \Phi_{d, \gm}(\bx, \bxi) \mu_\gm(\bxi)\,d\bxi. 
\eea
We will start with introducing our method for the one-dimensional cases, and then generalize it to two dimensions. 
The generalization to three dimensions is straightforward by following the same lines in two-dimensional scheme, and  we will omit it for brevity. 
Without loss of generality, we will derive the numerical scheme for parameter $\gm \in (\ap, 2]$.

\subsection{One-dimensional  cases}
\label{section2-1}

In one-dimensional ($d = 1$) cases, we let  domain $\Og = (a, b)$, and then region $\Upsilon = [0, L]$ and $\Upsilon_+^c = (L, \infty)$ with $L = b - a$.
For notational simplicity, we will  denote $\xi := \xi^{(1)}$ and $x := x^{(1)}$, and the function
\beas
\Phi_{1, \gm}(x,\, \xi) = \fl{u(x + \xi) - 2u(x) + u(x - \xi)}{\xi^\gamma}, \qquad \mbox{for \ $\gm \in (\ap, 2]$.}
\eeas
In the special case of $\gamma = 2$, function $\Phi_{1, 2}(x, \xi)$ represents the  central difference quotient for the classical Laplacian $\Dt$ at point $x$. 
 
Next, we will focus on approximating the integral ${\mathcal L}_{\Upsilon}^\ap u$.  
Define mesh size $h = (b - a)/N$ with a positive integer $N$, and denote grid points
\beas
x_i = a + ih, \qquad \xi_i = ih, \qquad \mbox{for} \ \ i = 0, 1, \ldots, N.
\eeas
It implies that $\xi_0 = 0$ and $\xi_{N} = L$. 
Assume the ansatz for function $\Phi_{1,\,\gamma}(x, \xi)$ as:
\begin{eqnarray}\label{interpolation}
\Phi_{1,\,\gamma}(x,\xi) = \sum _{k = 0}^{N} \Phi_{1,\,\gamma}(x,\xi_k)\,\varphi_k^p(\xi), \qquad \mbox{for} \ \ \xi \in [0, L],
\end{eqnarray}
where $\varphi_k^p(\xi)$, representing the $p$-th degree Lagrange polynomial, is the basis function at point $\xi_k$.  
The main novelty of our method is to interpolate function $\Phi_{1,\,\gamma}(x,\xi)$ with respect to $\xi$, instead of $x$.   
Note that at point $\xi = \xi_0$,  function $\Phi_{1,\,\gm}(x,\xi_0)$ is defined in a limit sense as  $\xi \rightarrow 0$. 
To see this, we will divide our discussion into two cases: 
\begin{itemize}
\item[(i).]  If the splitting parameter $\gm = 2$,   we obtain 
\begin{eqnarray}\label{gamma1}
\Phi_{1, 2}(x,\,\xi_0)
= \lim _{\xi \rightarrow 0} \frac{u(x-\xi)-2u(x)+u(x+\xi)}{\xi^2} = u''(x) \,\approx\, \Phi_{1,\,2}(x,\, \xi_1),\qquad 
\end{eqnarray}
via the central difference approximation to $u''(x)$. 
\item[(ii).] If $\gamma \in (\ap, 2)$,  there is
\begin{eqnarray}\label{gamma2}
\Phi_{1,\,\gamma}(x,\xi_0)
= \lim _{\xi \rightarrow 0} \frac{u(x-\xi)+u(x+\xi)-2u(x)}{\xi^2}\cdot \xi^{2-\gamma} = \Phi_{1,\,2}(x,\, \xi_0) \lim_{\xi \rightarrow 0}\xi^{2-\gamma}
=0.\quad 
\end{eqnarray}
Consequently,  the summation in (\ref{interpolation}) starts from $k = 1$ (instead of $k = 0$), if  $\gm \neq 2$.
\end{itemize}
Substituting the ansatz \eqref{interpolation} into  \eqref{L1} with $d = 1$, we obtain the approximation
\bea\label{FL3}
{\mathcal L}_{\Upsilon, h}^\ap u(x) = \sum_{k = 0}^{N} \Phi_{1,\,\gamma}(x, \xi_k) \int_0^L \varphi_k^p(\xi)\,\fl{d\xi}{\xi^{1+\ap-\gm}}, \qquad \mbox{for} \ \ x \in \Og.
\eea
For notational convenience, we will denote the {\it weight integral}
\begin{eqnarray}\label{wi-1D}
\og_k^p := \int_{0}^{L} \varphi_k^p(\xi)\,\fl{d\xi}{\xi^{1+\ap-\gm}},
\qquad \text{for}\ \  k = 0,1,2,\dots,N,
\end{eqnarray}
Since Lagrange polynomials are considered, the integral domain in \eqref{wi-1D} can be reduced from $[0, L]$ to the support of function $\varphi_k^p(\xi)$. 
Moreover, the weight integral $\og_k^p$ can be found analytically.
For convenience of the readers, we will provide the values of $\og_k^p$ for $p = 0, 1, 2$ in Appendix \ref{appendix1}. 
Substituting $x = x_i$ in \eqref{FL3} and noticing the definition of $\Phi_{1, \gm}(x, \xi)$, we obtain the approximation of ${\mathcal L}_{\Upsilon}^\ap u(x_i)$ as:
\begin{eqnarray}\label{L1h}
{\mathcal L}_{\Upsilon, h}^\ap u_i=\sum_{k = 1}^N\, \og_k^p\,\fl{u_{i-k} - 2u_i + u_{i+k}}{\xi_k^\gm}+\zeta\, \og_0^p\,\fl{u_{i-1} - 2u_i + u_{i+1}}{\xi_1^\gm}, \quad \ \mbox{for}  \ \ 1 \le i \le N-1,  
\end{eqnarray}
where we denote $u_i = u(x_i)$, and $\zeta = \lfloor{\gm}/{2}\rfloor$ with $\lfloor \cdot \rfloor$ representing the floor function.  
The second term in (\ref{L1h}) vanishes if $\gm \neq 2$, due to \eqref{gamma2}. 
On the other hand, the integral in \eqref{L2} at point $x = x_i$ becomes
\begin{eqnarray}\label{L2h}
{\mathcal L}_{\Upsilon_+^c}^\ap u_i =-\fl{2}{\ap L^\ap} u_i  + \int_{L}^\infty \big[g(x_i-\xi) + g(x_i + \xi)\big]\,\fl{d\xi}{\xi^{1+\ap}}, \quad \ \mbox{for}  \ \ 1 \le i \le N-1.
\end{eqnarray}

Combining  \eqref{L1h} and \eqref{L2h} and reorganizing the terms yield our numerical approximation to the one-dimensional integral fractional Laplacian $(-\Dt)^\fl{\ap}{2}$ with extended Dirichlet boundary conditions as:
\begin{eqnarray}\label{1D-scheme}
&&(-\Delta)^{\fl{\alpha}{2}}_hu_i
= -c_{1,\alpha}\bigg(a_0u_i + \sum_{k =1}^N a_k \big(u_{i+k}+u_{i-k}\big)\bigg)\nn\\
&&\hspace{1.6cm}= -c_{1,\alpha}\bigg(a_0u_i + \sum_{j = i+1}^{N-1} a_{j-i}\,u_j + \sum_{j = 1}^{i-1} a_{i - j} u_j + b_i\bigg), \qquad \mbox{for} \ \ 1 \le i \le N-1. \qquad\quad
\end{eqnarray}
The term $b_i$ comes from boundary conditions and is defined by
\beas
b_i = \sum_{j = N}^{N+i} a_{j-i}\, g(x_j)  + \sum_{j = i-N}^{0} a_{i-j}\,g(x_j) +  \int_{L}^\infty \big[g(x_i-\xi) + g(x_i + \xi)\big]\,\fl{d\xi}{\xi^{1+\ap}}.
\eeas
If homogeneous boundary conditions are considered, $b_i \equiv 0$ for $1 \le i \le N-1$.
The coefficients 
\beas
a_0 = -2\bigg(\sum_{j=1}^N a_j + \frac{1}{\alpha L^{\alpha}}\bigg), \qquad 
\mbox{with} \ \ a_j = \fl{1}{\xi_j^\gm}\left\{
\begin{array}{lcl}\displaystyle
\zeta\,\og_0^p + \og_1^p, &\ & \text{if} \ \  j = 1, \\ 
\displaystyle
\og_j^p,  & & \text{if}\ \  j \neq 0, 1.
\end{array}\right.
\eeas
The scheme in \eqref{1D-scheme} provides a general structure of our operator factorization method, and the dependence of basis function $\varphi_k^p$ is counted through the value of weight integral $\og^p_k$. 
In other words, if different basis functions are used, we only  need to update the value of $\og_k^p$ but keep the scheme and its computer implementation unchanged. 

Let ${\bf u} = \big(u_1, u_2,  \ldots,  u_{N-1}\big)^T$ and ${\bf b} = -c_{1, \ap}\big(b_1, b_2, \ldots,  b_{N-1}\big)^T$.
The discretization of the one-dimensional fractional Laplacian  in \eqref{1D-scheme} can be written into a matrix-vector form, i.e., $(-\Dt)^\fl{\ap}{2}_h{\bf u} = A^{(1)} {\bf u} + {\bf b}$,  where 
\beas
{{A}}^{(1)} = 
-c_{1,\alpha} \left(
\begin{array}{cccccc}
a_0 & a_1 & \ldots &  a_{N-3} & a_{N-2}  \\
a_1 & a_0 & a_1  &  \cdots & a_{N-3}  \\
\vdots & \ddots  & \ddots & \ddots & \vdots \\
a_{N-3} &  \ldots  & a_1 & a_0 & a_1 \\
a_{N-2} & a_{N-3}  & \ldots & a_1 & a_0
\end{array}
\right).
\eeas
It is clear that  $A^{(1)}$ is a symmetric Toeplitz matrix, and thus the product  ${A^{(1)} {\bf u}}$ can be computed efficiently by using the fast Fourier transform (FFT) with computational cost ${\mathcal O}\big((N-1) \log (N-1)\big)$.

\begin{remark}
Our method interpolating function $\Phi_{1, \gm}(x, \xi)$ with respect to $\xi$ is fundamentally different from that proposed in \cite{Huang2014}, although both of them use Lagrange basis functions $\varphi^p$ (for $p \in {\mathbb N}^0)$. 
The method in \cite{Huang2014}  has an $\ap$-dependent accuracy -- ${\mathcal O}(h^{2-\ap})$ for linear  basis functions $\varphi^1$,  while ${\mathcal O}(h^{3-\ap})$ for quadratic basis functions $\varphi^2$. 
In contrast to it, our method with $\gm = 2$ can achieve a uniform accuracy for any $\ap \in (0, 2)$. 
Specifically, the accuracy is ${\mathcal O}(h^{2})$  for both constant basis $\varphi^0$ and linear basis $\varphi^1$, while ${\mathcal O}(h^{4})$ for quadratic basis functions $\varphi^2$. 
More discussions and illustrations can be found in  Section \ref{section3-1}. 
\end{remark}
\subsection{Two-dimensional  cases}
\label{section2-2}

Our method uses Lagrange interpolations to the central difference quotient function and can be easily generalized to high dimensions with corresponding changes to ansatz  in (\ref{interpolation}). 
For convenience of the readers, we will show the two-dimensional scheme in this section, while the generalization to three dimensions will be omitted for brevity. 
In two-dimensional ($d = 2$) cases, we have 
\begin{eqnarray}\label{psi-2D}
\Phi_{2,\gamma}({\bf x},\textit{{\boldmath$\xi$}}) = 
\bigg(\sum_{\bm \in \varkappa_1}
u\big(x^{(1)}+(-1)^{m_1}\xi^{(1)}, \,x^{(2)}+(-1)^{m_2}\xi^{(2)}\big)-4 u({\bf x})
\bigg)\fl{1}{|\bxi|^\gm}.
\end{eqnarray}
Define mesh size $h = L/N$ with a positive integer $N$, and introduce grid points 
\beas
x_i^{(s)} = a_s + ih, \qquad \xi_i^{(s)} = ih, \qquad \mbox{for} \ \ i = 0, 1, \ldots, N, \quad s = 1, 2. 
\eeas
Denote point $\bxi_{kl}  = \big(\xi_k^{(1)}, \,\xi_l^{(2)}\big)$, for $1 \le k, l \le N$. 
We assume the two-dimensional ansatz: 
\begin{eqnarray}\label{interpolation-2D}
\Phi_{2, \gamma}({\bf x},\,\textit{{\boldmath$\xi$}}) = 
\sum_{k = 0}^N\sum_{l = 0}^N \Phi_{2, \gamma}({\bf x},\,\bxi_{kl})\,\varphi^p_{k}\big(\xi^{(1)}\big) \varphi^p_{l}\big(\xi^{(2)}\big), \qquad \mbox{for} \ \  \bxi \in \Upsilon.
\end{eqnarray}
Similar to the one-dimensional cases, function $\Phi_{2, \gm}(\bx, \bxi_{00})$ is defined in a limit sense as $\bxi \to (0, 0)$. 
Specifically, we approximate it as 
\beas
\Phi_{2, 2}(\bx, \bxi_{00}) &=& \lim_{\bxi \to (0, 0)} \Phi_{2, 2}(\bx, \bxi)\nn\\
&\approx& \Phi_{2, 2}(\bx, \bxi_{10}) + \Phi_{2, 2}(\bx, \bxi_{01}) -\Phi_{2, 2}(\bx, \bxi_{11}), \qquad \mbox{for} \ \ \gm = 2. 
\eeas
For $\gm \in (\ap, 2)$, we get
\beas
\Phi_{2, \gm}(\bx, \bxi_{00})  = \lim_{\bxi \to {\bf 0}}\Phi_{2, \gm}(\bx,\,\bxi) = \Phi_{2, 2}(\bx,\,\bxi) \lim_{\bxi \to {\bf 0}}|\bxi|^{2-\gm} = 0, \qquad \mbox{for} \ \ \gm \in (\ap, 2). 
\eeas
Substituting ansatz \eqref{interpolation-2D} into ${\mathcal L}_{\Upsilon}^\ap u$ with $d = 2$, we obtain 
\bea\label{2D-scheme0}
{\mathcal L}_{\Upsilon, h}^\ap u(\bx) = 
\sum_{k  = 0}^N \sum_{l = 0}^N \Phi_{2, \gm}(\bx,\,\bxi_{kl})\underbrace{\int_0^L \int_0^L \varphi_k^p\big(\xi^{(1)}\big)\varphi_l^p\big(\xi^{(2)}\big)\fl{d\bxi}{|\bxi|^{2+\ap-\gm}}}_{\og_{kl}^p}, \qquad \mbox{for} \ \ \bx \in \Og. 
\eea
The weight integral $\og_{kl}^p$ is actually on the intersection region of the supports of two basis function $\varphi_k^p$ and $\varphi_l^p$, instead of $[0, L] \times [0, L]$. 
Different from one-dimensional cases, it is challenging to obtain the analytical results of these integrals, and thus we will use numerical quadrature rules to compute them. 
Noticing the approximation of $\Phi_{2, \gm}(\bx, \bxi_{00})$, we can further formulate (\ref{2D-scheme0}) as: 
 \beas
&&{\mathcal L}_{\Upsilon, h}^\ap u(\bx) 
= \mathop{\sum_{k=1}^{N} \sum_{l=1}^{N}}_{k+l\neq 2} \og_{kl}^p\,\Phi_{2, \gm}(\bx, \bxi_{kl}) + \sum_{k = 2}^N \Big(\og_{0k}^p\,\Phi_{2, \gm}(\bx, \bxi_{0k}) + \og_{k0}^p\,\Phi_{2, \gm}(\bx, \bxi_{k0}) \Big)\nn\\
&&\hspace{1.7cm} + \big(\og_{01} + \zeta\og_{00}\big)\Phi_{2, \gm}(\bx, \bxi_{10}) + \big(\og_{01} + \zeta\og_{00}\big)\Phi_{2, \gm}(\bx, \bxi_{10}) + \big(\og_{11} - \zeta\og_{00}\big)\Phi_{2, \gm}(\bx, \bxi_{11}).\qquad
\eeas

Without loss of generality, we assume that $N_1 = N$, and choose $N_2$  as the smaller integer such that  $a_2 + N_2h \geq b_2$. 
Denote $\bx_{ij} = \big(x_i^{(1)}, x_j^{(2)}\big)$ and $u_{ij} = u(\bx_{ij})$. 
Substituting $\bx = \bx_{ij}$  into ${\mathcal L}_{\Upsilon, h}^\ap u(\bx)$ and  ${\mathcal L}_{\Upsilon_+^c}^\ap u(\bx)$, noticing the definition of $\Phi_{2, \gm}(\bx, \bxi)$, we then obtain the numerical approximation the two-dimensional  fractional Laplacian $(-\Dt)^\fl{\ap}{2}$ with Dirichlet boundary conditions as follows: 
 \bea\label{2D-scheme}
&&{\mathcal L}_{\Upsilon, h}^\ap u_{ij}
= -c_{2,\alpha}\bigg[a_{00}\,u_{ij} + \sum_{k = 0}^{i-1}\bigg(\sum_{\substack{l = 0 \\ k+l\neq 0}}^{j-1} a_{kl}u_{(i-k)(j-l)} + 
\sum_{\substack{l = 1}}^{N_2 - (j+1)} a_{kl}u_{(i-k)(j+l)}\bigg)\nn\\
&&\hspace{3.9cm}+ \sum_{k = 1}^{N_1 - (i+1)}\bigg(\sum_{l = 0}^{j-1} a_{kl}u_{(i+k)(j-l)} + 
\sum_{\substack{l = 1}}^{N_2 - (j+1)} a_{kl}u_{(i+k)(j+l)}\bigg) + b_{ij}\bigg],\qquad\qquad 
\eea
for $1 \le i \le N_1-1$ and $1 \le j \le N_2-1$, 
where the coefficients 
\beas
&&a_{kl} = \fl{1}{|\bxi_{kl}|^{\gm}}
\left\{\begin{array}{ll}
\displaystyle 2\og_{kl}^p +\zeta \og_{00}^p,  \quad \  & \mbox{if \  $(k, l) = (0, 1)$\ or\ $(1, 0)$}, \\
\displaystyle \og_{kl}^p -\zeta \og_{00}^p,  \quad \  & \mbox{if \  $k = l = 1$}, \\
\displaystyle 2\og_{kl}^p,   &  \mbox{if} \ \ k = 0 \ \ \mbox{or} \ \, l = 0, \ \ \mbox{and} \ \  k + l \ge 2, \\ 
\og_{kl}^p,    &  \mbox{otherwise},
\end{array} \right.\qquad\qquad\\
&&a_{00}
= -2\sum_{k=1}^{N} \big(a_{k0} + a_{0k}\big)  -4\sum_{i=1}^{N}\sum_{j=1}^{N} a_{ij}
-4 \int_{\Upsilon_+^c}\fl{d\bxi}{|\bxi|^{2+\ap}}. \nn
\eeas
The term $b_{ij}$ is defined as
\beas
&&b_{ij} = \sum_{l = -N}^N \bigg(\sum_{k = i}^N a_{k|l|}\,g\big(\bx_{(i-k)(j+l)}\big) + \sum_{k = N_1-i}^{N} a_{k|l|}\,g\big(\bx_{(i+k)(j+l)}\big)\bigg) + \sum_{k = 0}^{i-1}\bigg(\sum_{l = j}^N a_{kl}\,g\big(\bx_{(i-k)(j-l)}\big)\qquad \nn\\
&&\hspace{0.8cm}+\sum_{l = N_2-j}^N a_{kl}\,g\big(\bx_{(i-k)(j+l)}\big)\bigg)+ \sum_{k = 1}^{N_1-(i+1)}\bigg(\sum_{l = j}^N a_{kl}\,g\big(\bx_{(i-k)(j-l)}\big) +\sum_{l = N_2-j}^N a_{kl}\,g\big(\bx_{(i-k)(j+l)}\big)\bigg) \nn\\
&&\hspace{0.8cm}+ \sum_{{\bf m}\in \varkappa_1}\int_{\Upsilon_+^c}g\big(\bx_{ij} + (-1)^{\bf m}\circ\bxi\big)\fl{d\bxi}{|\bxi|^{2+\ap}},
\eeas
for $1 \le i \le N_1-1$ and $1 \le j \le N_2-1$. 

We can further write the scheme in (\ref{2D-scheme}) into matrix-vector form. 
 For $1\le j \le N_2-1$, denote the vector ${\bf u}_{j}^{(1)} = \big(u_{1j}, u_{2j},\ldots,u_{(N_2-1)j}\big)$, and let the block vector $
{\bf u} = \big({\bf u}_{1}^{(1)}, {\bf u}_{2}^{(1)},\,\ldots,{\bf u}_{N_1-1}^{(1)}\big)^T$, and the block vector ${\bf b}$ is defined in the same manner as ${\bf u}$ with entries $-c_{2, \ap}b_{ij}$. 
Then the matrix-vector form of the  scheme (\ref{2D-scheme}) is given by 
$(-\Delta)^{\fl{\alpha}{2}}_{h}{\bf u} = A^{(2)}{\bf u} + {\bf b}$, where $A^{(2)}$ is a block-Toeplitz--Toeplitz-block  matrix defined as
\beas\label{A-2D}
{ {A}}^{(2)}= 
\left(
\begin{array}{cccccc}
A_{0} & A_{1} & \ldots &  A_{N_2-3} & A_{N_2-2}   \\
A_{1} & A_{0}  & A_{1}   &  \cdots & A_{N_2-3}   \\
\vdots & \ddots  & \ddots & \ddots & \vdots \\
A_{N_2-3}  &  \ldots  & A_{1}  & A_{0}  & A_{1}  \\
A_{N_2-2}  & A_{N_2-3}  & \ldots & A_{1}  & A_{0} 
\end{array}
\right)
\eeas
with each block $A_j$ a symmetricToeplitz matrix
\beas\label{Aj}
{\rm {A}}_{j} = 
\left(
\begin{array}{cccccc}
a_{0j} & a_{1j} & \ldots &  a_{(N_1-3)j} & a_{(N_1-2)j}  \\
a_{1j} & a_{0j} & a_{1j}  &  \cdots & a_{(N_1-3)j}  \\
\vdots & \ddots  & \ddots & \ddots & \vdots \\
a_{(N_1-3)j} &  \ldots  & a_{1j} & a_{0j} & a_{1j} \\
a_{(N_1-2)j} & a_{(N_1-3)j}  & \ldots & a_{1j} & a_{0j}
\end{array}
\right), \qquad\mbox{for \ $j = 0, 1, \ldots, N_2-2$. }
\eeas 
It is easy to see that  $A^{(2)}$ is a full positive definite matrix,  which usually requires huge memory and computational costs for computing matrix-vector products. 
However, the Toeplitz structure of $A^{(2)}$  enables us to develop fast algorithms via fast Fourier transform and thus efficiently compute matrix-vector multiplication. 

\begin{remark}[Extension to other nonlocal operators]
Our method provides a general framework of operator factorization, and can be easily generalized to solve the nonlocal operators of the form: 
\bea\label{Remark00}
{\mathcal L}u(\bx) = \bar{c}_{d, \ap}\int_{{\mathbb R}^d}\fl{u(\bx) - u({\bf y})}{|\bx - {\bf y}|^{d+\ap}} K(|\bx - {\bf y}|) d{\bf y},
\eea
where $K(r)$ denotes a  kernel function, e.g., $K(r) = \exp(-\lambda r)$ in the tempered fractional Laplacian \cite{Duo-TFL2019}. 
In this case, we can formulate the operator \eqref{Remark00} as 
\beas\label{Remark01}
{\mathcal L}u(\bx) = -\bar{c}_{d, \ap}\int_{{\mathbb R}_+^d}\Phi_{d, \gm}(\bx, \bxi) \fl{K(|\bxi|)}{|\bxi|^{d+\ap - \gm}}d\bxi
\eeas
with $\Phi_{d, \gm}(\bx, \xi)$ defined in (\ref{eq0}). 
We then follow the same procedure as for \eqref{FL2} to approximate it. 
Note that the extra kernel function changes the function $\mu_\gm(\bxi)$ and thus affects the values of  $\og^p$, but the scheme structure  remains the same as that of fractional Laplacian $(-\Dt)^\fl{\ap}{2}$.
\end{remark}

\section{Numerical experiments}
\label{section3}
\setcounter{equation}{0}

In this section, we test the performance of our method in discretizing the fractional Laplacian $(-\Dt)^\fl{\ap}{2}$ and in solving fractional Poisson equations. 
Numerical results from  constant basis ($\varphi^0$), linear basis ($\varphi^1$), and quadratic basis ($\varphi^2$) functions are compared and discussed under different conditions of function $u$ and power $\ap$. 
Moreover,  nonhomogeneous Dirichlet boundary conditions will be considered, which have been rarely studied in the literature  \cite{Burkardt0020, Wu0021, Acosta2019}. 
Unless otherwise stated, we will choose splitting parameter $\gamma = 2$ in our simulations. 
More discussion and comparisons of different splitting parameters can be found in Example 3.1.3.

\subsection{Estimation of the fractional Laplacian}
\label{section3-1}

First, we study the accuracy of our method in approximating the fractional Laplacian $(-\Dt)^{\fl{\ap}{2}}$  on a bounded domain with extended  Dirichlet  (homogeneous or nonhomogeneous) boundary conditions. 
Here, we consider the one-dimensional cases with $\Og = (-1,1)$, and thus Dirichlet boundary conditions are imposed on $\Og^c = (-\infty, -1]\cup[1, \infty)$. 
Denote the error function  as
\beas
e_{\Dt}(x) = (-\Dt)^{\fl{\ap}{2}}_h u(x) - (-\Dt)^{\fl{\ap}{2}} u(x), \qquad \mbox{for} \ \ x \in \Og,
\eeas
where $(-\Dt)^{\fl{\ap}{2}}_h$ represents a numerical approximation of the fractional Laplacian $(-\Dt)^\fl{\ap}{2}$. 
Then numerical accuracy under different smoothness conditions of $u$ will be studied in  the following examples.

\bigskip 
\noindent{\bf Example 3.1.1 \big(Less smooth function $u$).\ } 
Here, we are interested in understanding the minimum conditions for our method to be consistent. 
To this end, we consider a compact support function $u(x) = (1-x^2)_+^{1+\lfloor \ap\rfloor}$. 
That is,  function $u(x)  = (1-x^2)_+$ for $\ap < 1$, while $u(x) = (1-x^2)_+^2$ for $1 \le \ap < 2$.
Table \ref{Table3-1-1} shows numerical errors $\|e_\Dt\|_\infty$ and convergence rates for different basis function $\varphi^p$ and power $\ap$.
\begin{table}[htb!]
        \centering
        \begin{tabular}{|c||c|c|c|c|c|c|}
        \hline
        $h$ & 1/16 & 1/32 & 1/64 & 1/128 & 1/256 & 1/512 \\
        \hline
        \multicolumn{7}{|c|}{$\ap = 0.5$}\\
        \hline
        \multirow{2}{*}{$\varphi^0$}& 7.5879e-3 & 5.3220e-3 & 3.7499e-3 & 2.6475e-3 & 1.8707e-3 &1.3224e-3\\
        \cline{2-7}
        & \mbox{c.r.}& 0.5117 & 0.5051 & 0.5023 & 0.5010& 0.5005\\
        \hline
        \multirow{2}{*}{$\varphi^1$}&7.8596e-3 &  5.4986e-3  & 3.8696e-3 &  2.7304e-3  & 1.9287e-3  & 1.3632e-3\\
        \cline{2-7}
        & \mbox{c.r.}&0.5154  & 0.5069  & 0.5031 &  0.5014 &  0.5007\\
        \hline
        \multirow{2}{*}{$\varphi^2$}&1.3338e-2  & 9.3477e-3 &  6.5851e-3 &  4.6488e-3 &  3.2848e-3  & 2.3219e-3\\
        \cline{2-7}
        & \mbox{c.r.}&0.5129  & 0.5054 &  0.5024 &  0.5011  & 0.5005\\
        \hline
        \multicolumn{7}{|c|}{$\ap = 1$}\\
        \hline
        \multirow{2}{*}{$\varphi^0$}&8.1722e-4 &  3.8342e-4  & 1.8911e-4 &  9.4360e-5 &  4.7189e-5 &  2.3604e-5\\
        \cline{2-7}
        & \mbox{c.r.}&1.0918 &  1.0197  & 1.0029  &0.9997  & 0.9994\\
        \hline
        \multirow{2}{*}{$\varphi^1$}&8.1722e-4  & 3.8342e-4 &  1.8911e-4&   9.4360e-5   &4.7189e-5  & 2.3604e-5\\
        \cline{2-7}
        & \mbox{c.r.}&1.0918 &  1.0197 &  1.0029 &  0.9997  & 0.9994\\
        \hline
        \multirow{2}{*}{$\varphi^2$}&4.7698e-3 & 2.3572e-3 &  1.1717e-3  & 5.8413e-4 &  2.9164e-4  & 1.4571e-4\\
        \cline{2-7}
        & \mbox{c.r.}&1.0169&  1.0085 &  1.0042  & 1.0021&   1.0011\\
        \hline
        \multicolumn{7}{|c|}{$\ap = 1.7$}\\
        \hline
        \multirow{2}{*}{$\varphi^0$}&3.6356e-3  & 1.8041e-3  & 1.2777e-3 &  1.0195e-3 &  8.3276e-4  & 6.8083e-4\\
        \cline{2-7}
        & \mbox{c.r.}&1.0109 &  0.4977  & 0.3257  & 0.2919 &  0.2906\\
        \hline
        \multirow{2}{*}{$\varphi^1$}&2.5288e-3 &  5.4873e-4 &  5.8948e-4 &  5.0041e-4  & 4.0137e-4 &  3.2097e-4\\
        \cline{2-7}
        & \mbox{c.r.}&2.2043 & -0.1033 & 0.2363 &  0.3182 &  0.3225\\
        \hline
        \multirow{2}{*}{$\varphi^2$}&9.9878e-2 &  7.8950e-2 &  6.3253e-2  & 5.1024e-2 &  4.1302e-2 &  3.3489e-2\\
        \cline{2-7}
        & \mbox{c.r.}&0.3392 &  0.3198  & 0.3010 &  0.3050 &  0.3025\\
        \hline
        \end{tabular}
        \caption{Numerical errors $\|e_\Dt\|_\infty$ and convergence rates (c.r.) in approximating function  $(-\Dt)^\fl{\ap}{2}u$ with $u = (1-x^2)_+^{1+\lfloor\ap \rfloor}$ and basis function $\varphi^p$ (for $p = 0, 1, 2$).} \label{Table3-1-1}
    \end{table}

It shows that as mesh size $h$ reduces, numerical errors decrease with a rate depending on power $\ap$. 
Specifically, all three basis functions $\varphi^p$ (for $p = 0, 1, 2$) lead to  the same convergence rate -- ${\mathcal O}(h^{1-\ap})$ for $\ap < 1$, while ${\mathcal O}(h^{2-\ap})$ for $\ap \ge 1$. 
Numerical errors are maximized around the boundary points $x = \pm 1$. 
In this case, the function $u$ has limited smoothness at two points $x = \pm 1$, which creates a bottleneck for numerical accuracy (see similar observations of finite difference methods in \cite{Duo2018, Duo-FDM2019}).
\begin{figure}[htb!]
\centerline{(a) \includegraphics[height=5.660cm,width=7.560cm]{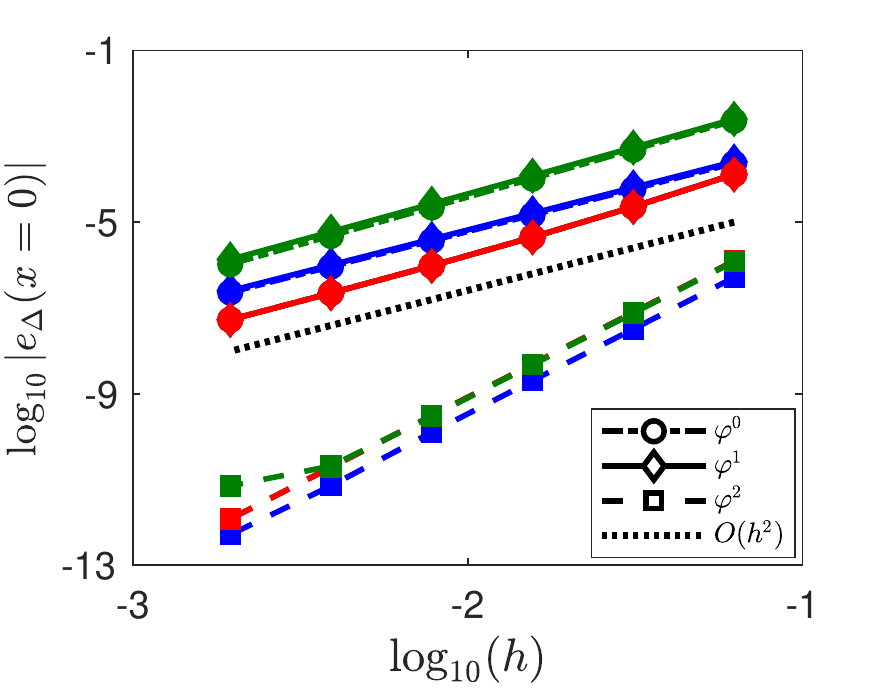}
\hspace{-5mm}
(b) \includegraphics[height=5.660cm,width=7.560cm]{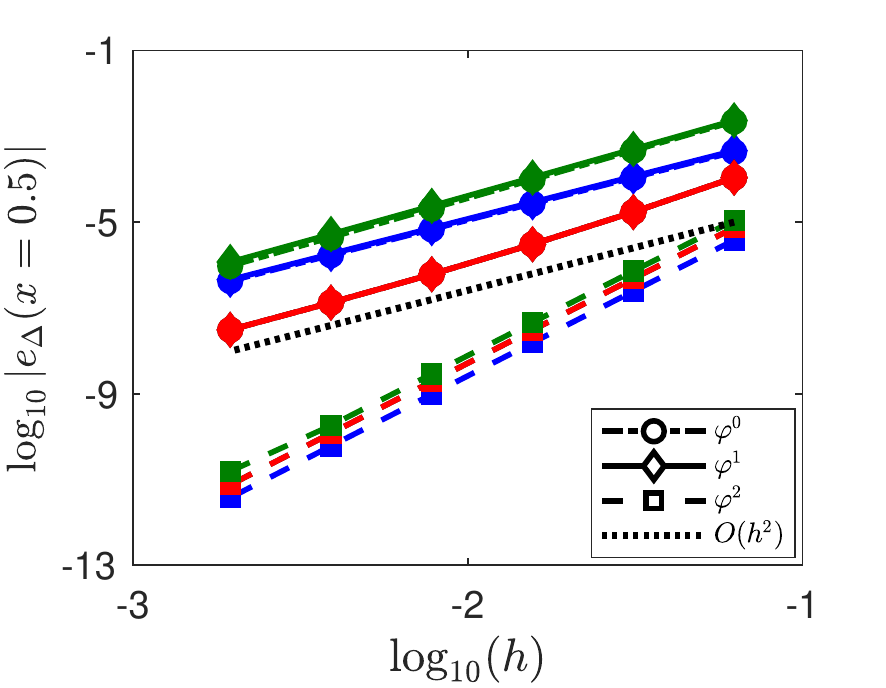}}
\caption{Numerical error $|e_\Dt(x)|$ at point $x = 0$ (left) and $x = 0.5$ (right) in approximating function $(-\Dt)^\fl{\ap}{2}u$ with $u(x) = (1-x^2)_+^{1+\lfloor\ap\rfloor}$, where $\ap = 0.6$ (blue), $1$ (red), and $1.5$ (green).}\label{Figure3-1-1}
\end{figure}
But, our method can achieve the optimal accuracy --  ${\mathcal O}(h^2)$ for constant basis $\varphi^0$ and linear basis $\varphi^1$, while ${\mathcal O}(h^4)$ for quadratic basis $\varphi^2$ -- at inner points.  
To illustrate it,  Figure \ref{Figure3-1-1} presents numerical errors at point $x = 0$ and $x = 0.5$, where an order line (i.e., dotted line) is included for easy comparison.  
We find that numerical errors from constant basis and linear basis are almost the same for different $\ap$, and they are  identical if $\ap = 1$. 
Numerical errors of quadratic basis $\varphi^2$ at inner points (e.g. $x = 0, 0.5$) are much smaller than those from constant and linear bases, suggesting the benefits of using high-degree basis functions. 
Moreover, comparing Figure \ref{Figure3-1-1} (a) and (b) shows that numerical error $|e_\Dt(x)|$ increases as $x$ approaches the domain boundary, consistent with our earlier observations of maximum errors around boundary points. 

The above results and our extensive studies suggest that  the minimum consistency condition of our method in approximating the Dirichlet fractional Laplacian $(-\Dt)^\fl{\ap}{2}$ depends on power $\ap$. 
For $m \in {\mathbb N}$ and $l \in (0, 1]$, we denote 
\beas
C^{m,\,l}(\bar{\Og}) = \bigg\{u \in C^{m}(\bar{\Og}) \ \bigg| \ \sup_{\substack{x,\, y \in \bar{\Og} \\ x \neq y}}\fl{|u^{(k)}(x) - u^{(k)}(y)|}{|x - y|^l} < \infty, \ \ \mbox{for} \ \, k \in{\mathbb N} \ \, \mbox{and} \ \, k\le m\bigg\}.
\eeas
Then we can summarize the consistency results as follows. 
\begin{remark}[{\bf Consistency conditions}]
\label{remark3-1}
Let $(-\Dt)^{\fl{\ap}{2}}_{h}$ be a numerical approximation to the Dirichlet fractional Laplacian $(-\Dt)^\fl{\ap}{2}$  with  small mesh size $h$. 
Our method with basis function $\varphi^p$ (for $p = 0, 1,$ or $2$) is consistent if 
function $u \in C^{\lfloor\ap\rfloor,\, \ap-\lfloor\ap\rfloor+\veps}(\bar{\Og})$ with small number $\veps > 0$. 
Moreover, their local truncation errors satisfy 
\beas 
\|e_\Dt\|_\infty \le Ch^\veps,\qquad \mbox{for}  \ \ \ap \in (0, 2)
\eeas 
with $C$ a positive constant independent of $h$.
\end{remark}
Our extensive studies show that the consistency result in Remark \ref{remark3-1} is independent of splitting parameter $\gm$, i.e., it holds for any $\gm \in (\ap, 2]$.

\bigskip 
\noindent{\bf Example 3.1.2 \big(Smooth function $u$\big). \ } 
In this example, we continue our study to test the accuracy of our method for functions that are smoother around boundary points.  
We are interested in understanding the optimal accuracy of different basis functions and also the minimum conditions to achieve such an optimal accuracy. 

To this end, let's first consider function $u = (1-x^2)_+^{2.1+\ap}$. 
Table \ref{Table3-1-2} presents numerical errors $\|e_\Dt\|_\infty$ and convergence rates for different basis functions $\varphi^p$.
\begin{table}[h!]
        \centering
        \begin{tabular}{|c||c|c|c|c|c|c|}
        \hline
        $h$ & 1/16 & 1/32 & 1/64 & 1/128 & 1/256 & 1/512 \\
               \hline
        \multicolumn{7}{|c|}{$\ap = 0.5$}\\
        \hline
        \multirow{2}{*}{$\varphi^0$}&9.7624e-5 &  2.8827e-5&   7.6872e-6 &  1.9687e-6  & 4.9667e-7 &  1.2457e-7\\
        \cline{2-7}
        & \mbox{c.r.}&1.7598&   1.9069 &  1.9653&   1.9868 &  1.9953\\
        \hline
        \multirow{2}{*}{$\varphi^1$}&2.1391e-4 &  5.9663e-5 &  1.5540e-5 &  3.9426e-6  & 9.9077e-7&   2.4817e-7\\
        \cline{2-7}
        & \mbox{c.r.}&1.8421 &  1.9409  & 1.9787 &  1.9925 &  1.9972\\
        \hline
        \multirow{2}{*}{$\varphi^2$}&1.0716e-4 &  2.4168e-5 &  5.5431e-6 &  1.2823e-6  & 2.9789e-7&  
        6.9347e-8 \\
        \cline{2-7}
        & \mbox{c.r.}&2.1486 &  2.1243  & 2.1120  & 2.1059 &  2.1029\\
        \hline
        \multicolumn{7}{|c|}{$\ap = 1$}\\
        \hline
        \multirow{2}{*}{$\varphi^0$}&5.9137e-4  & 7.5126e-5  & 9.4842e-6  & 2.0487e-6 &  6.4898e-7  & 1.7163e-7\\
        \cline{2-7}
        & \mbox{c.r.}&2.9767 &  2.9857 &  2.2109&   1.6584  & 1.9189\\
        \hline
        \multirow{2}{*}{$\varphi^1$}&5.9137e-4  & 7.5126e-5  & 9.4842e-6 &  2.0487e-6&   6.4898e-7 &  1.7163e-7\\
        \cline{2-7}
        & \mbox{c.r.}&2.9767 &  2.9857 &  2.2109  & 1.6584  & 1.9189\\
        \hline
        \multirow{2}{*}{$\varphi^2$}&2.4438e-4&   5.4962e-5 &  1.2583e-5  & 2.9079e-6  & 6.7516e-7  & 1.5713e-7\\
        \cline{2-7}
        & \mbox{c.r.}&2.1526 &  2.1270 &  2.1134 &  2.1067 &  2.1033\\
        \hline
        \multicolumn{7}{|c|}{$\ap = 1.7$}\\
        \hline
        \multirow{2}{*}{$\varphi^0$}&1.4385e-2  & 5.3211e-3 &  1.5476e-3  & 4.0205e-4  & 9.9477e-5  & 2.4066e-5\\
        \cline{2-7}
        & \mbox{c.r.}&1.4348 &  1.7817 &  1.9446  & 2.0150  & 2.0474\\
        \hline
        \multirow{2}{*}{$\varphi^1$}&1.5206e-2  & 5.5501e-3  & 1.6318e-3 &  4.2867e-4  & 1.0729e-4  & 2.6267e-5\\
        \cline{2-7}
        & \mbox{c.r.}&1.4540  & 1.7660 &  1.9285  & 1.9983  & 2.0302\\
        \hline
        \multirow{2}{*}{$\varphi^2$}&8.0506e-4 &  1.0070e-4  & 1.6011e-5 &  3.1025e-6  & 6.6674e-7  & 1.4996e-7\\
        \cline{2-7}
        & \mbox{c.r.}&2.9991  & 2.6529 &  2.3675 &  2.2182  & 2.1526\\
        \hline
        \end{tabular}
        \caption{Numerical errors $\|e_\Dt\|_\infty$ and convergence rate (c.r.) in approximating function $(-\Dt)^\fl{\ap}{2}u$ with $u = (1-x^2)_+^{2.1+\ap}$ and basis function $\varphi^p$ (for $p = 0, 1, 2$).}\label{Table3-1-2}
        \label{tab2}
    \end{table}
Compared to Table \ref{Table3-1-1}, numerical errors in this case are much smaller since function $u$ is smoother at boundary. 
But,  maximum errors are still found around two boundary points $x = \pm 1$. 
Both constant basis $\varphi^0$ and linear basis $\varphi^1$ have accuracy rate of ${\mathcal O}(h^2)$, while the quadratic basis $\varphi^2$ leads to a higher rate, i.e., ${\mathcal O}(h^{2.1})$.
Moreover, numerical errors from quadratic basis are much smaller than those from constant and linear bases. 

\begin{figure}[htb!]
\centerline{(a)\includegraphics[height=5.660cm,width=7.560cm]{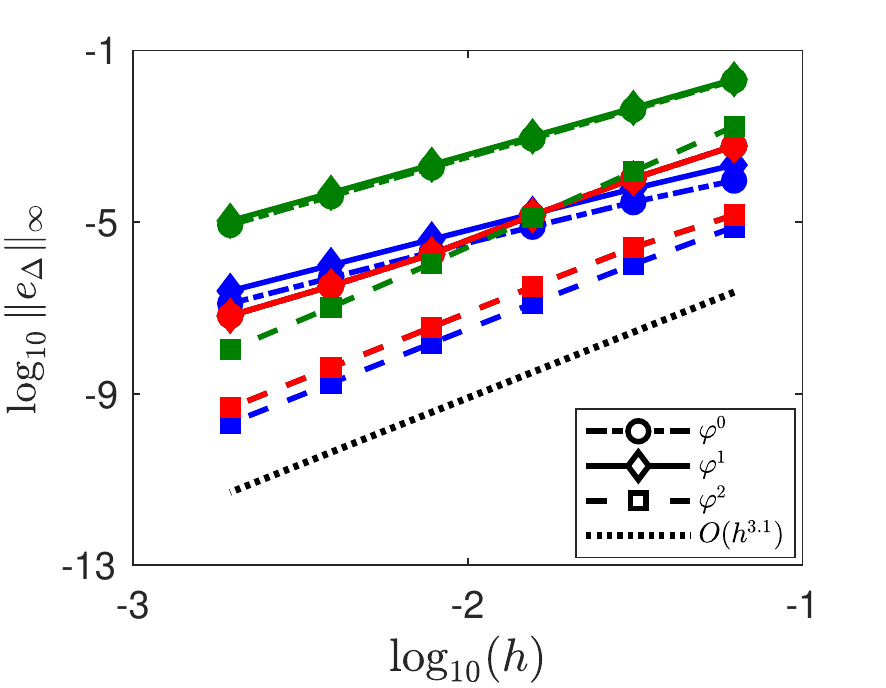}\hspace{-5mm}
(b)\includegraphics[height=5.660cm,width=7.560cm]{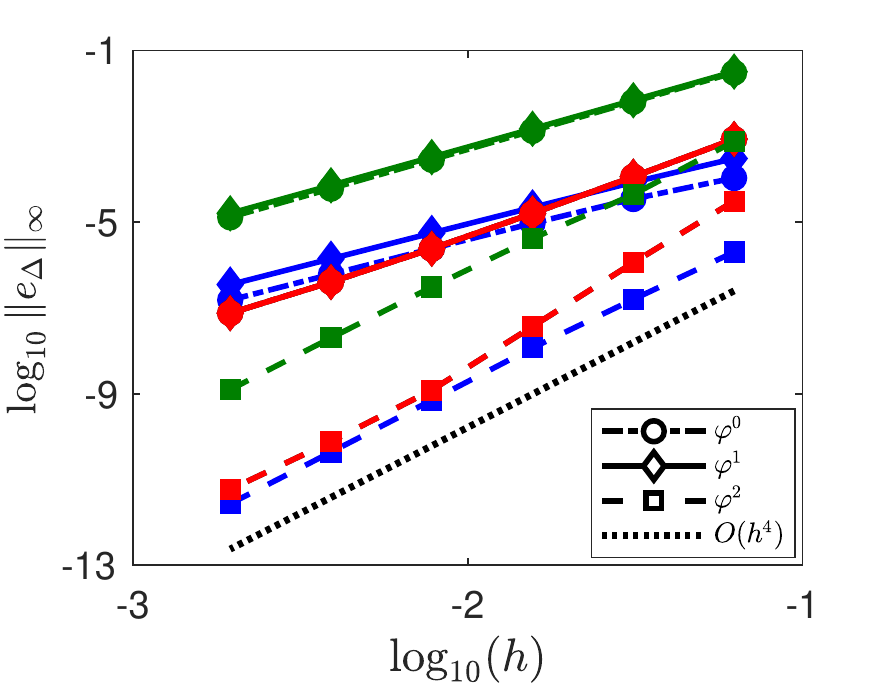}}
\caption{Numerical errors $\|e_\Dt\|_\infty$ in approximating function $(-\Dt)^\fl{\ap}{2}u$, where $\ap = 0.6$ (blue), $1$ (red), and $1.5$ (right). (a) $u = (1-x^2)_+^{3.1+\ap}$; (b)  $u = (1-x^2)_+^{4.1+\ap}$.}\label{Figure3-1-2}
\end{figure}

To further our understanding, we increase the smoothness of functions around boundary points.  
Figure \ref{Figure3-1-2} shows numerical errors for functions $u(x) = (1-x^2)_+^{3.1+\ap}$ and $u = (1-x^2)_+^{4.1+\ap}$.
For both functions, the constant basis $\varphi^0$ and linear basis $\varphi^1$ remain the second-order accuracy, suggesting the best accuracy of these two basis functions is ${\mathcal O}(h^2)$. 
In contrast, the accuracy of quadratic basis $\varphi^2$ increases as function $u$ becomes smoother, and the best accuracy is  ${\mathcal O}(h^4)$ (see Figure \ref{Figure3-1-2} (b)). 
Figure \ref{Figure3-1-2} additionally suggests that the smaller the power $\ap$, the less the numerical errors. 
From the above observations and our extensive studies, we summarize the optimal accuracy of our methods in the following remark. 

\begin{remark}[{\bf Optimal accuracy}]
\label{Remark3-2}
Let $(-\Dt)^{\fl{\ap}{2}}_{h, 2}$ be a numerical approximation to the Dirichlet fractional Laplacian $(-\Dt)^\fl{\ap}{2}$  with  splitting parameter $\gm = 2$ and small mesh size $h$.  
\begin{itemize}
\item[(i) ]  If function $u \in C^{2+\lfloor\ap\rfloor, \,\ap-\lfloor\ap\rfloor+\veps}(\bar{\Og})$, the error function  from basis function $\varphi^p$ satisfies
\beas
\|e_\Dt\|_\infty \le C h^2, \qquad \mbox{for}\ \ \ap \in (0, 2), 
\eeas
for $p = 0, 1, 2$. 
\item[(ii) ] If function $u \in C^{4+\lfloor\ap\rfloor, \,\ap-\lfloor\ap\rfloor+\veps}(\bar{\Og})$, the error function from quadratic basis $\varphi^2$  satisfies
\beas
\|e_\Dt\|_\infty \le C h^4, \qquad \mbox{for}\ \ \ap \in (0, 2).
\eeas
\end{itemize}
Here, $C$ is a positive constant independent of mesh size $h$, and $\veps > 0$ is a small number. 
\end{remark}
Our studies show  that the optimal accuracy can be achieved only when splitting parameter $\gm = 2$. 
In contrast, other splitting parameters $\gm \in (\ap, 2)$ lead to larger numerical errors and lower accuracy rates. 
More comparisons and discussions of different splitting parameters can be found in next example.

\bigskip 
\noindent{\bf Example 3.1.3 \big(Infinitely smooth function $u$\big).  } 
In this example, we consider infinitely smooth function $u = 1/(1+x^2)$ for $x \in {\mathbb R}$,  and  approximate  $(-\Dt)^\fl{\ap}{2}u(x)$ for $x \in (-1, 1)$. 
It is equivalent to approximate the fractional Laplacian with extended nonhomogeneous boundary conditions $u = 1/(1+x^2)$ for $x \in (-1, 1)^c$.   
In addition to numerical accuracy, we will also study the effects of splitting parameter $\gm$. 

Table \ref{Table3-1-3} shows that quadratic basis functions yield significantly smaller numerical errors than  constant and linear bases. 
This is one main advantage of our method -- when function $u$ is smooth enough, our method enables us to  increase the accuracy by using high-degree basis functions (e.g., $\varphi^2$). 
\begin{table}[h!]
        \centering
        \begin{tabular}{|c||c|c|c|c|c|}
        \hline
        $h$ & 1/16 & 1/32 & 1/64 & 1/128 & 1/256  \\
               \hline
        \multicolumn{6}{|c|}{$\ap = 0.5$}\\
        \hline
        \multirow{2}{*}{$\varphi^0$}&2.0023e-5  & 5.4222e-6  & 1.3919e-6  & 3.5121e-7 &  8.8084e-8 \\
        \cline{2-6}
        & \mbox{c.r.}&1.8847 &  1.9618  & 1.9866 &  1.9954  \\
        \hline
        \multirow{2}{*}{$\varphi^1$}&3.0000e-5 &  8.3090e-6  & 2.1718e-6 &  5.5132e-7 &  1.3857e-7 \\
        \cline{2-6}
        & \mbox{c.r.}& 1.8522 &  1.9357 &  1.9780 &  1.9923  \\
        \hline
        \multirow{2}{*}{$\varphi^2$}&1.7384e-7 &  1.2148e-8  & 7.8796e-10&   4.9569e-11 &  2.6986e-12 \\
        \cline{2-6}
        & \mbox{c.r.}&3.8391 &  3.9464 &  3.9906 &  4.1991  \\
        \hline
        \multicolumn{6}{|c|}{$\ap = 1$}\\
        \hline
        \multirow{2}{*}{$\varphi^0$}&1.1056e-4  & 1.7994e-5  & 3.2863e-6 &  6.6985e-7 &  1.4849e-7  \\
        \cline{2-6}
        & \mbox{c.r.}&2.6193 &  2.4530  & 2.2945 &  2.1735  \\
        \hline
        \multirow{2}{*}{$\varphi^1$}&1.1056e-4 &  1.7994e-5 &  3.2863e-6  & 6.6985e-7  & 1.4849e-7\\
        \cline{2-6}
        & \mbox{c.r.}&2.6193  & 2.4530   &2.2945  & 2.1735 \\
        \hline
        \multirow{2}{*}{$\varphi^2$}&8.0637e-7 &  2.5951e-8  & 8.3813e-10  & 2.7735e-11  & 9.3070e-13  \\
        \cline{2-6}
        & \mbox{c.r.}&4.9576  & 4.9525  & 4.9174 &  4.8972 \\
        \hline
        \multicolumn{6}{|c|}{$\ap = 1.7$}\\
        \hline
        \multirow{2}{*}{$\varphi^0$}&2.2284e-3 &  4.6838e-4  & 9.8834e-5  & 2.0988e-5 &  4.4908e-6 \\
        \cline{2-6}
        & \mbox{c.r.}&2.2503 &  2.2446  & 2.2355 &  2.2245  \\
        \hline
        \multirow{2}{*}{$\varphi^1$}&2.3784e-3 &  5.1283e-4  & 1.1135e-4 &  2.4403e-5 &  5.4026e-6\\
        \cline{2-6}
        & \mbox{c.r.}&2.2135 &  2.2033 &  2.1900 &  2.1753  \\
        \hline
        \multirow{2}{*}{$\varphi^2$}&3.1706e-5 &  1.6865e-6  & 8.9337e-8  & 4.7546e-9 &  2.6046e-10  \\
        \cline{2-6}
        & \mbox{c.r.}&4.2326 &  4.2386 &  4.2319  & 4.1902  \\
        \hline
        \end{tabular}
        \caption{Numerical errors $\|e_\Dt\|_\infty$ and convergence rate (c.r.) in approximating function $(-\Dt)^\fl{\ap}{2}u$ on (-1, 1) with $u(x) = 1/(1+x^2)$ and basis function $\varphi^p$ (for $p = 0, 1, 2$).}\label{Table3-1-3}
    \end{table}
In terms of computer implementation, different basis functions only change the predefined values of $\og_k^p$, but do not affect the structures of main programs. 
Additionally, numerical errors in this example are maximized at point $x = 0$,  different from those in Examples 3.1.1 and 3.1.2.  
Note that even though both $\varphi^0$ and $\varphi^1$ have the same accuracy rate ${\mathcal O}(h^2)$,  numerical errors from constant basis are generally smaller  if $\ap \neq 1$. 

\begin{figure}[htb!]
\centerline{
\includegraphics[height = 5.66cm, width = 7.56cm]{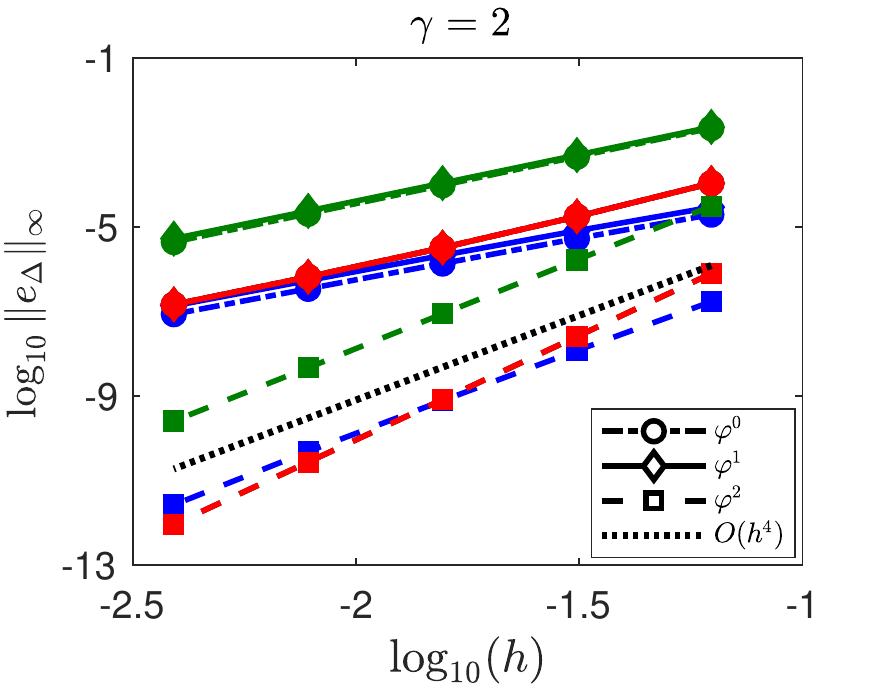}
\includegraphics[height = 5.66cm, width = 7.56cm]{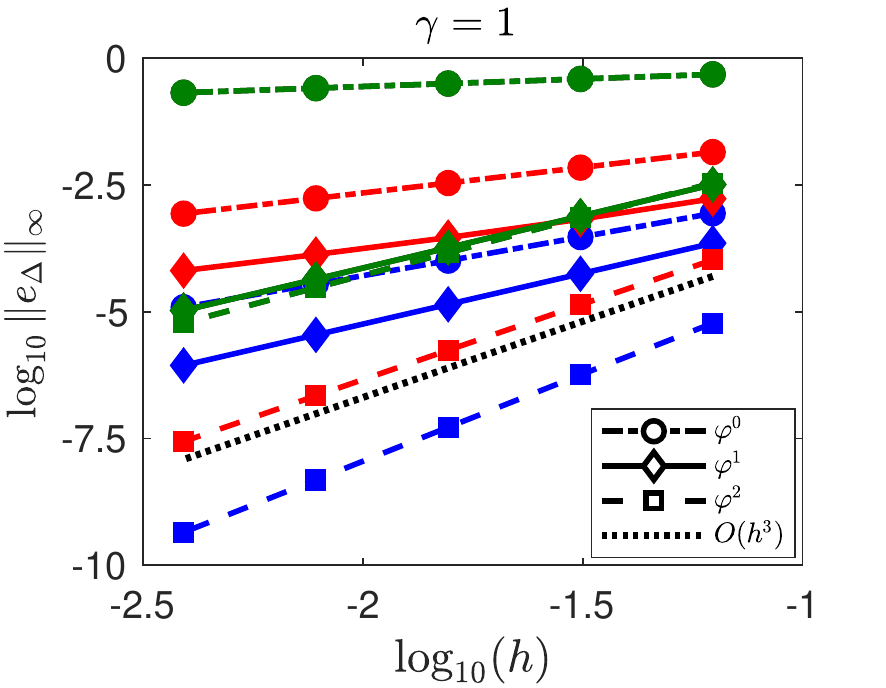}}
\vspace{3mm}
\centerline{
\includegraphics[height = 5.66cm, width = 7.56cm]{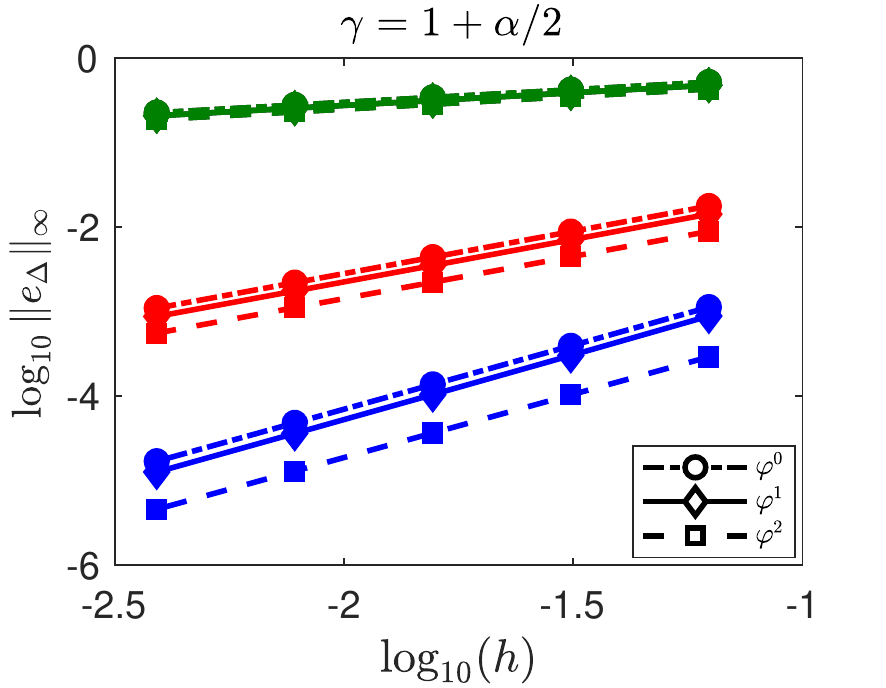}
\includegraphics[height = 5.66cm, width = 7.56cm]{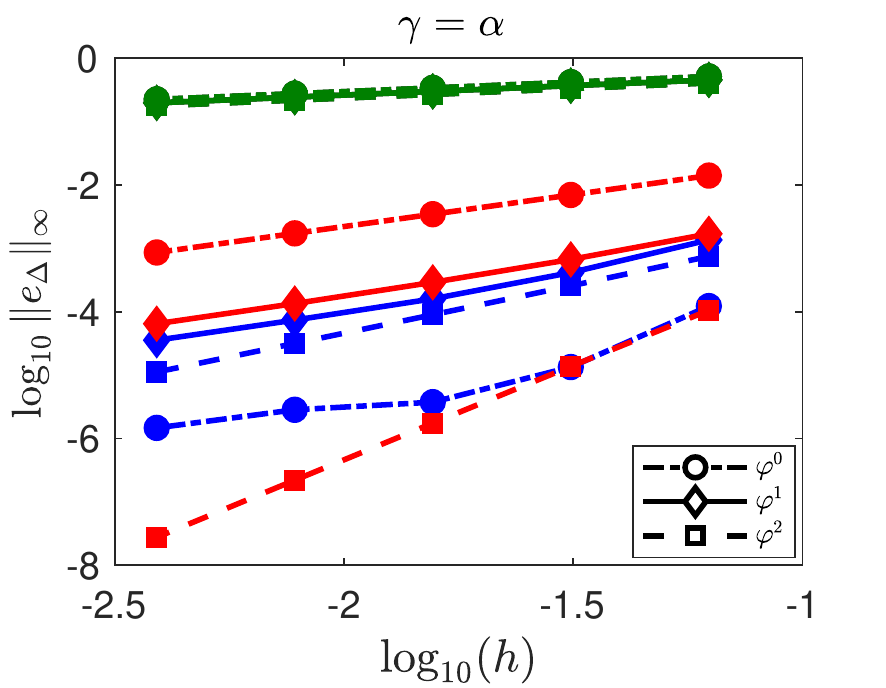}}
\caption{Effects of splitting parameter $\gm$  in approximating function $(-\Dt)^\fl{\ap}{2}u$ on $(-1, 1)$ with $u(x) = 1/(1+x^2)$, where $\ap = 0.6$ (blue), $1$ (red), and $1.5$ (right).}\label{Figure3-1-3}
\end{figure}

In Figure \ref{Figure3-1-3}, we study the effects of splitting parameter $\gm$ on numerical accuracy for different basis function $\varphi^p$, where $\gm = 2$, $1$, $1+\ap/2$, and $\ap$ are considered. 
It shows that the splitting parameter $\gm = 2$ leads to the smallest numerical errors for each basis function. 
Moreover, when $\gm = 2$ our method can achieve the optimal convergence rates, i.e., ${\mathcal O}(h^2)$ for  $\varphi^0$ and $\varphi^1$, while ${\mathcal O}(h^4)$ for  $\varphi^2$. 
In contrast to $\gm = 2$, an $\ap$-dependent  rate is observed if  splitting parameter $\gm \neq 2$. 
For instance, when $\gm = 1+\ap/2$ we find that the accuracy is ${\mathcal O}(h^{2-\ap})$ for constant basis $\varphi^0$ and linear basis $\varphi^1$, while ${\mathcal O}(h^{3-\ap})$ for  quadratic basis $\varphi^2$.  
Generally, the accuracy rate for $\gm \in (0, 2)$ depends on power $\ap$ and is lower than that from $\gm = 2$.

\subsection{Solution of Poisson problems}
\label{section3-2}

In this section, we further explore the performance of our method in solving the fractional Poisson equation and its generalization. 
Denote the error function in solution as
\beas
e_u(\bx) = u_h(\bx) - u(\bx),  \qquad \mbox{for} \ \ \bx\in\Og,
\eeas
where $u_h$ and $u$ represent numerical and exact solutions, respectively. 
As seen in Example 3.1.3,  the splitting parameter $\gm = 2$ outperforms  other parameter $\gm \in (\ap, 2)$, and thus it will be always used in the following examples. 

\bigskip
\noindent{\bf Example 3.2.1 \big(1D fractional Poisson equation\big).  } 
Let the one-dimensional domain $\Og  = (-1, 1)$.  
First, we consider the benchmark fractional Poisson problem, i.e., choosing $f(x)  = 1$ in \eqref{Poisson} and $g(x) \equiv 0$ in \eqref{BC}.   
Its  exact solution is given by 
$$u(x) = -\fl{1}{\Gamma(1+\ap)}(1-x^2)^{\fl{\ap}{2}}_+, \qquad \mbox{for} \ \ x\in {\mathbb R},$$
i.e., $u \in C^{0, \fl{\ap}{2}}(\bar{\Og})$.  
Table \ref{Table3-2-1} presents numerical errors $\|e_u\|_\infty$ and convergence rates for different basis function $\varphi^p$ and power $\ap$.

\begin{figure}[htb!]
\centerline{\includegraphics[height=5.660cm,width=7.560cm]{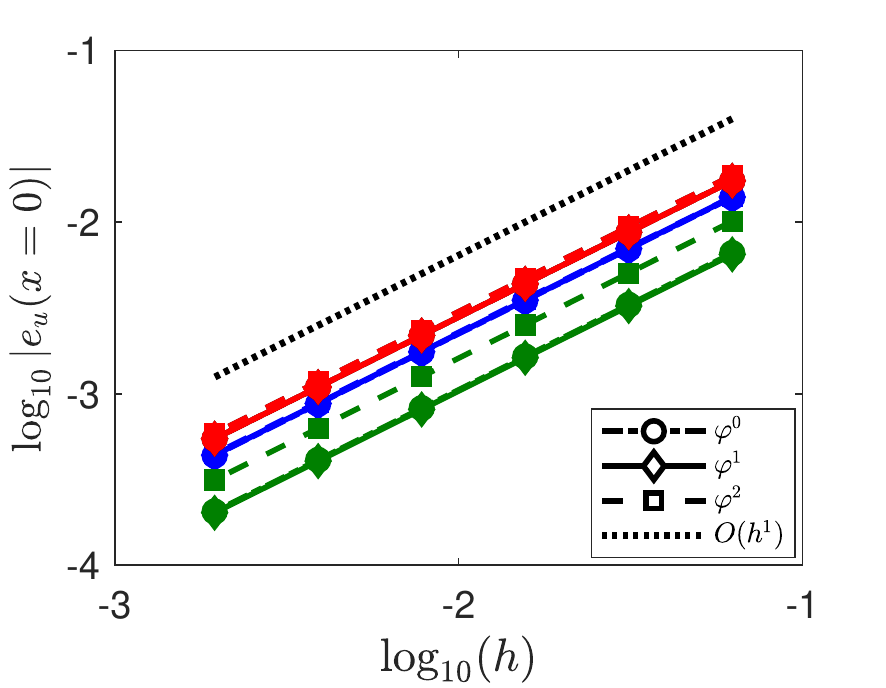}}
\caption{Numerical errors $|e_u(x)|$ at point $x = 0$ in solving the 1D fractional Poisson problem \eqref{Poisson}--\eqref{BC} with $f(x) = 1$ and $g(x) = 0$, where $\ap = 0.6$ (blue), $1$ (red), or $1.5$ (green).}\label{Figure3-2-1}
\end{figure}

It shows that our method has an accuracy rate of ${\mathcal O}(h^\fl{\ap}{2})$ for any $\ap \in (0, 2)$, independent of  basis function $\varphi^p$. 
In this case, the regularity of solution at boundary creates a bottleneck for numerical methods, and consequently using high-degree basis function does not necessarily benefit the accuracy. 
Note that the same convergence rate is observed for finite difference methods \cite{Duo2018, Duo-FDM2019, Duo-TFL2019}. 
Generally, the larger the power $\ap$, the smoother the solution at boundary, the smaller the numerical errors. 
Similar to the observations in Section \ref{section3-1}, the constant basis $\varphi^0$ and linear basis $\varphi^1$ have the same numerical errors if $\ap = 1$, {since their differentiation matrices for the fractional Laplacian $(-\Dt)^\fl{1}{2}$ are identical.} 
 \begin{table}[htb!]
        \centering
        \begin{tabular}{|c||c|c|c|c|c|c|}
        \hline
        $h$ & 1/16 & 1/32 & 1/64 & 1/128 & 1/256 & 1/512 \\
        \hline
        \multicolumn{7}{|c|}{$\ap = 0.6$}\\
        \hline
        \multirow{2}{*}{$\varphi^0$}&7.4494e-2 &  5.9980e-2 &  4.8507e-2 &  3.9314e-2 &  3.1898e-2 &  2.5895e-2\\
        \cline{2-7}
        & \mbox{c.r.}& 0.3126&   0.3063 &  0.3031 &  0.3016 &  0.3008 \\
        \hline
        \multirow{2}{*}{$\varphi^1$}&7.5493e-2 &  6.0790e-2  & 4.9164e-2 &  3.9847e-2 &  3.2331e-2 &  2.6247e-2\\
        \cline{2-7}
        & \mbox{c.r.}& 0.3125 &  0.3062  & 0.3031 &  0.3016  & 0.3008\\
        \hline
        \multirow{2}{*}{$\varphi^2$}&8.4532e-2&   6.8106e-2   &5.5102e-2  & 4.4671e-2&   3.6249e-2 &  2.9429e-2\\
        \cline{2-7}
        & \mbox{c.r.}& 0.3117 &  0.3057 &  0.3028  & 0.3014 &  0.3007\\
        \hline
        \multicolumn{7}{|c|}{$\ap = 1$}\\
        \hline
        \multirow{2}{*}{$\varphi^0$}&4.9166e-2 &  3.4508e-2  & 2.4310e-2 &  1.7158e-2 &  1.2121e-2  & 8.5671e-3\\
        \cline{2-7}
        & \mbox{c.r.}& 0.5107  & 0.5054 &  0.5027  & 0.5013 &  0.5007\\
        \hline
        \multirow{2}{*}{$\varphi^1$}&4.9166e-2  & 3.4508e-2  & 2.4310e-2  & 1.7158e-2  & 1.2121e-2 &  8.5671e-3\\
        \cline{2-7}
        & \mbox{c.r.}&0.5107 &  0.5054 &  0.5027 &  0.5013  & 0.5007\\
        \hline
        \multirow{2}{*}{$\varphi^2$}&5.7935e-2 &  4.0695e-2  & 2.8682e-2  & 2.0248e-2  & 1.4306e-2 &  1.0112e-2\\
        \cline{2-7}
        & \mbox{c.r.}& 0.5096  & 0.5047   & 0.5023  & 0.5012 &  0.5006\\
        \hline
        \multicolumn{7}{|c|}{$\ap = 1.5$}\\
        \hline
        \multirow{2}{*}{$\varphi^0$}&1.6161e-2 &  9.5429e-3 &  5.6545e-3  & 3.3563e-3 &  1.9939e-3 &  1.1851e-3\\
        \cline{2-7}
        & \mbox{c.r.}& 0.7600 &  0.7550 &  0.7525   & 0.7513 &  0.7506 \\
        \hline
        \multirow{2}{*}{$\varphi^1$}&1.5976e-2 &  9.4344e-3  & 5.5905e-3 &  3.3184e-3  & 1.9714e-3  & 1.1717e-3\\
        \cline{2-7}
        & \mbox{c.r.}& 0.7599  & 0.7550 &  0.7525  & 0.7513 &  0.7506\\
        \hline
        \multirow{2}{*}{$\varphi^2$}&2.2627e-2 &  1.3365e-2  & 7.9205e-3 &  4.7018e-3 &  2.7934e-3 &  1.6603e-3\\
        \cline{2-7}
        & \mbox{c.r.}& 0.7596 &  0.7548 &  0.7524 &  0.7512  & 0.7506\\
        \hline
        \end{tabular}
        \caption{Numerical errors $\|e_u\|_\infty$  and convergence rate (c.r.) in solving the 1D Poisson problem on $\Og = (-1, 1)$,  where $f(x) = 1$ in \eqref{Poisson} and $g(x) = 0$ in \eqref{BC}.}\label{Table3-2-1}
    \end{table}        
On the other hand,  the maximum numerical errors in Table \ref{Table3-2-1} are found around two boundary points $x = \pm 1$. 
But,  numerical errors at points far away from the domain boundary are much smaller, and moreover an accuracy rate of ${\mathcal O}(h)$ is observed at these points; see Figure. \ref{Figure3-2-1} for errors at $x = 0$.

\begin{figure}[htb!]
\centerline{(a)\includegraphics[height=5.660cm,width=7.560cm]{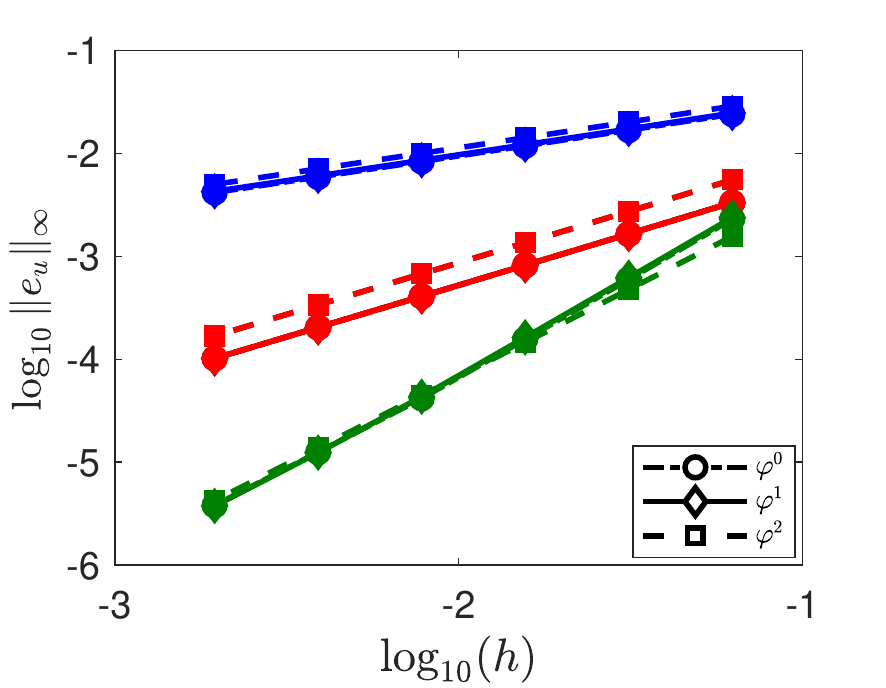}
(b)\includegraphics[height=5.660cm,width=7.560cm]{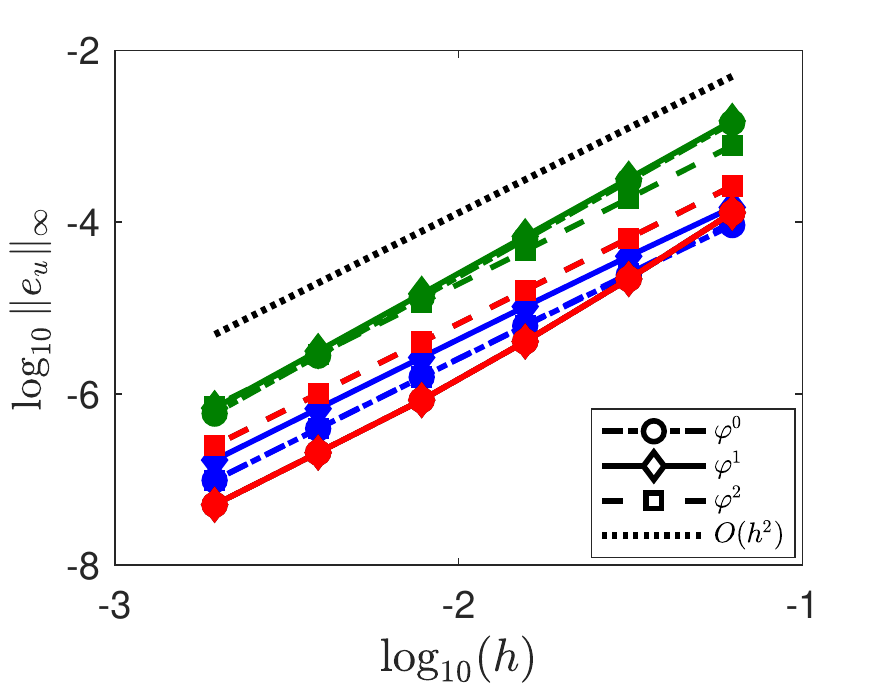}}
\centerline{(c)\includegraphics[height=5.660cm,width=7.560cm]{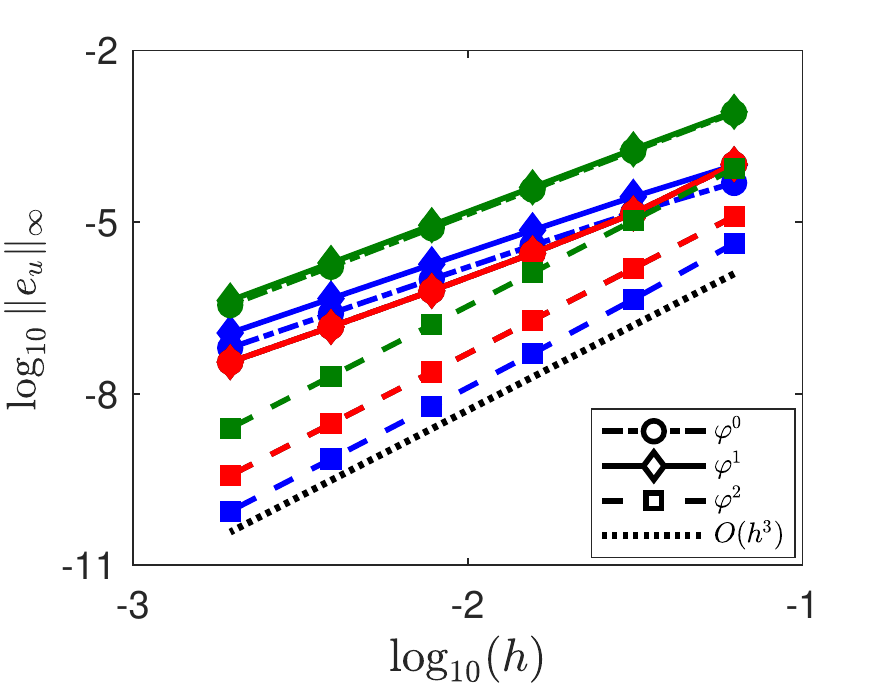}
(d)\includegraphics[height=5.660cm,width=7.560cm]{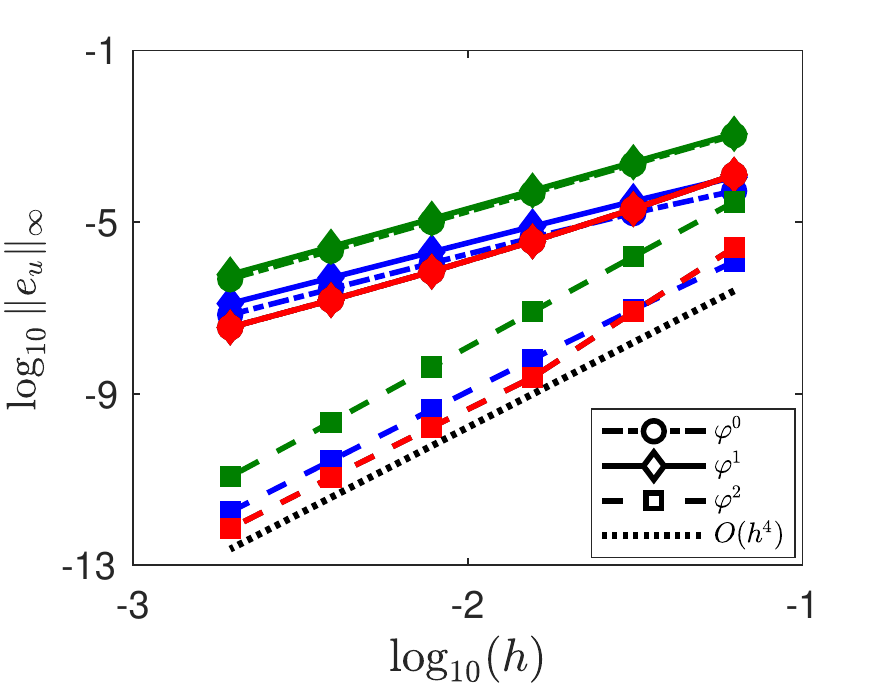}}
\caption{Numerical errors $\|e_u\|_\infty$ in solving the 1D Poisson problem \eqref{Poisson}--\eqref{BC} with $g(x) = 0$ and $f(x)$ in (\ref{fun0}), where the exact solution is  $u(x) = (1-x^2)_+^s$. 
From (a) to (d): $s = \ap, \, 2, \, 3,\, 4$, where $\ap = 0.6$ (blue), $1$ (red), or $1.5$ (green). }\label{Figure3-2-2}
\end{figure}

Next, we generalize our study  and consider  the Poisson problem \eqref{Poisson}--\eqref{BC} with extended homogeneous boundary conditions $g(x) \equiv 0$ and 
\bea\label{fun0}
f(x) = \fl{2^\ap\Gamma(\fl{\ap+1}{2})\Gamma(s+1)}{\sqrt{\pi}\,\Gamma(s+1-\fl{\ap}{2})}\,_2F_1\Big(\fl{\ap+1}{2}, -s + \fl{\ap}{2}; \, \fl{1}{2}; \, x^2\Big), \qquad \mbox{for} \ \ x \in (-1, 1),
\eea
for $s > 0$, where $\,_2F_1$ denotes the hypergeometric Gauss function.  
The  exact solution in this case is given by $u(x) = (1-x^2)_+^s$. 
It is clear that (\ref{fun0}) is a general case of the  benchmark problem (e.g., $s = \ap/2$), and the regularity of solutions can be controlled by parameter $s$ -- the larger the value of $s$, the smoother the solution at boundary. 
Figure \ref{Figure3-2-2} shows the numerical errors for different parameters $s = \ap, 2, 3,$ and $4$. 
From it, we find that  
\begin{enumerate}\itemsep -2pt
\item[i).] If $s \le 2$, all three basis functions lead to an $s$-dependent accuracy, i.e.,  $\|e_u\|_\infty \sim {\mathcal O}(h^s)$.  
In this case, the regularity of solution plays a dominant role in determining the accuracy of numerical methods (see similar observations for finite difference method in \cite[Table3]{Duo2018}.)
\item[ii).] If $s > 2$, the solution $u \in C^{2}(\bar{\Og})$. 
In this case, both constant basis $\varphi^0$ and linear basis $\varphi^1$ functions achieve the optimal  accuracy of $\|e_u\|_\infty \sim {\mathcal O}(h^2)$. 
Moreover, their numerical errors are almost the same (see Figure \ref{Figure3-2-2} (b)--(d)). 
In contrast, the quadratic basis $\varphi^2$ function can achieve a higher accuracy rate ${\mathcal O}(h^s)$ for $2 < s <  4$ (see Figure \ref{Figure3-2-2} (c)--(d)).
\item[iii).] If $s \ge 4$, the quadratic basis $\varphi^2$ has much smaller numerical errors than those two basis functions and reach the optimal accuracy rate $\|e_u\|_\infty \sim {\mathcal O}(h^4)$.
\end{enumerate}
The above studies suggest when choosing basis functions, one should take the regularity of solution into account so as to achieve the best performance of our method.

\bigskip
\noindent{\bf Example 3.2.2 \big(1D tempered fractional Poisson equation\big).  } 
Remark 2.2 shows that our method can be easily generalized  to study other nonlocal operators of form in (\ref{Remark00}). 
To illustrate this,  we will apply our method to study the Poisson equation with tempered fractional Laplacian, i.e. choosing kernel function $K(|x - y|) = e^{-\lambda|x-y|}$ in \eqref{Remark00}.  
The tempered fractional Laplacian is used to study the coexistence of normal and anomalous diffusion in many complex systems, where $\lambda \ge 0$ is a model parameter that mediates these two  diffusion   (see \cite{Duo-TFL2019} and references therein for more discussion). 
In our example, the right hand side function $f(x)$ is chosen such that the exact solution  is given by $u(x) = (1-x^2)^2_+$ for $x \in {\mathbb R}$.

\begin{figure}[htb!]
\centerline{(a)\includegraphics[height=5.660cm,width=7.560cm]{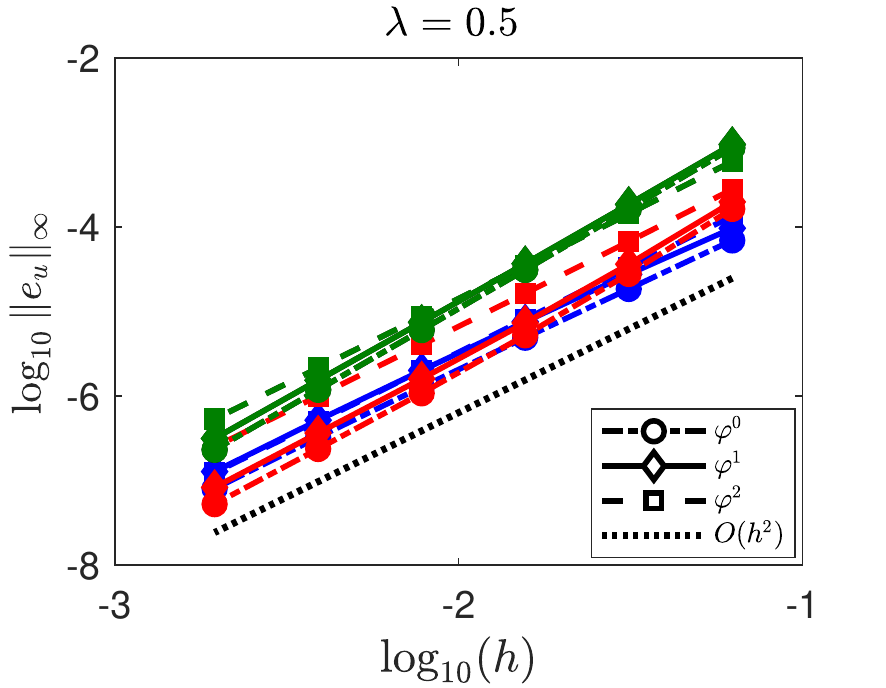}
(b)\includegraphics[height=5.660cm,width=7.560cm]{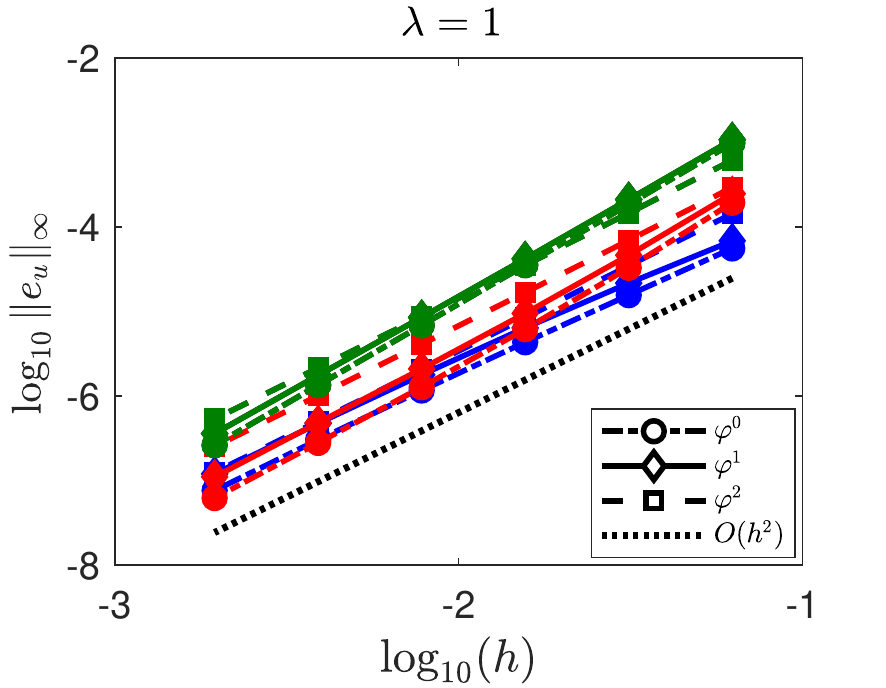}}
\caption{Numerical errors $\|e_u\|_\infty$ in solving the 1D tempered fractional Poisson problem with exact solution $u(x) = (1-x^2)_+^2$, where $\ap = 0.6$ (blue), $1$ (red), or $1.5$ (green). }\label{Figure3-2-0}
\end{figure}

Figure \ref{Figure3-2-0} shows numerical errors for $\lambda = 0.5$ and $\lambda = 1$. 
Note that when $\lambda =0$, it reduces to the fractional Poisson problem studied in Figure \ref{Figure3-2-2} (b). 
It shows that the uniform second-order accuracy is achieved for both $\lambda = 0.5$ and $\lambda = 1$, same rate as in Figure \ref{Figure3-2-2} (b) for $\lambda = 0$. 
Moreover, the numerical errors are insensitive to parameter $\lambda$.  
The computer implementation of this problem is the same as that for the fractional Poisson problem in Example 3.2.1, except the values of $\og_k^p$. 
Due to  the tempered term $e^{-\lambda|x-y|}$, the integrals of $\og^p_k$ in (\ref{wi-1D}) cannot be  analytically obtained, but they can be pre-calculated with numerical quadrature rules. 

\bigskip
\noindent{\bf Example 3.2.3 \big(2D fractional Poisson equation\big). } 
In this example, we consider the two-dimensional Poisson problem  \eqref{Poisson}--\eqref{BC} on a square domain $\Og = (-1, 1)^2$ with extended nonhomogeneous boundary conditions:
\beas\label{fun1}
g(\bx) = e^{-|\bx|^2}, \qquad \ \mbox{for} \ \ \bx \in \Og^c.
\eeas
The right hand function in \eqref{Poisson} is chosen as 
\beas
f(\bx) = 2^\ap\Gamma\Big(1+\fl{\ap}{2}\Big)\,_1F_1\Big(1+\fl{\ap}{2};\,1;\, -|\bx|^2\Big), \qquad\,\mbox{for} \ \ \bx \in \Og, \eeas
where $\,_1F_1$ denotes the confluent hypergeometric function. 
Compared to homogeneous cases, nonhomogeneous boundary conditions require extra computation of integrals involving nonzero $g(\bx)$, e.g., the second term in ${\mathcal L}_{\Upsilon_+^c}^\ap$, which can be accurately computed with traditional quadrature rules. 
Since the discretization matrix is of multi-level Toeplitz structure, the resulting linear 
system can be solved with fast algorithms based on the fast Fourier transforms. 

\begin{figure}[htb!]
\centerline{\includegraphics[height=4.860cm,width=5.360cm]{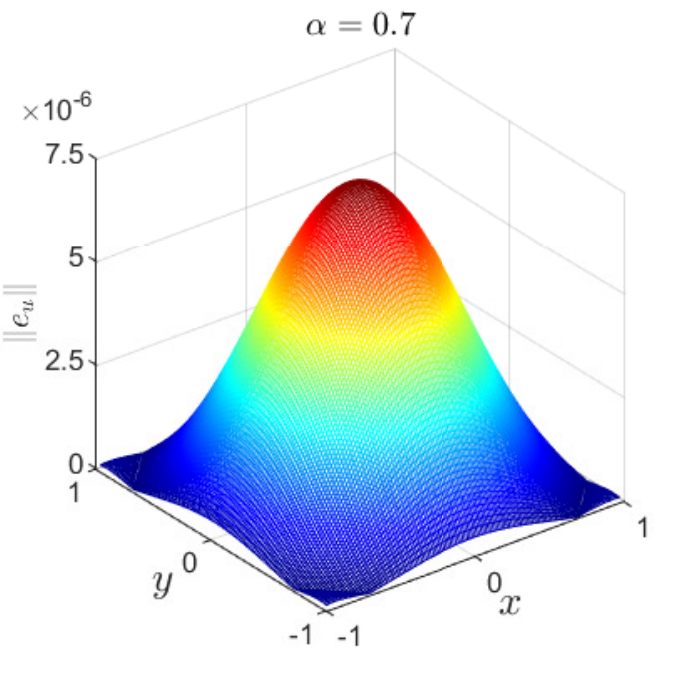}\hspace{5mm}
\includegraphics[height=4.860cm,width=5.360cm]{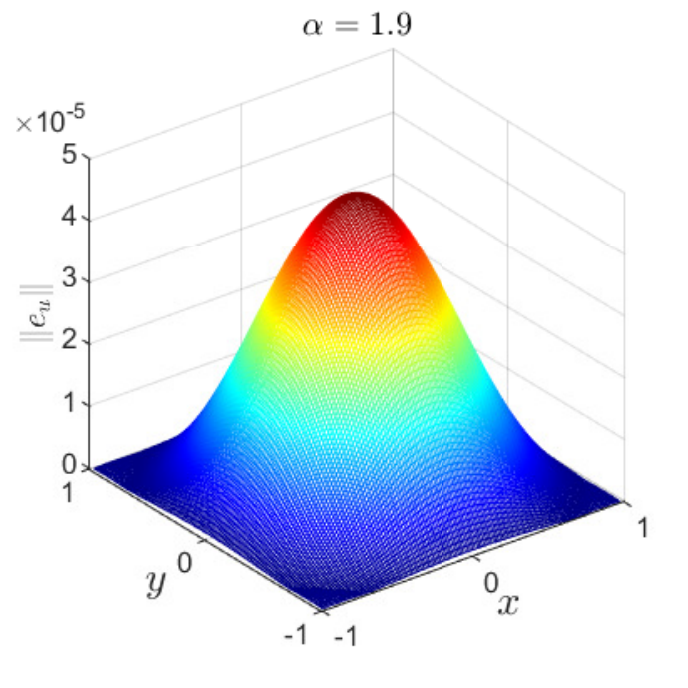}}
\caption{Numerical errors in the solution of the 2D fractional Poisson problems with basis function $\varphi^1$ and mesh size $h = 1/64$.}\label{Figure3-2-3}
\end{figure}
Figure \ref{Figure3-2-3} shows numerical errors in solution for $\ap = 0.7$ and $1.9$ with linear basis function $\varphi^1$ and  mesh size $h = 1/64$. 
It is clear that numerical errors reduce as $|\bx|$ increases, and  the maximum errors are obtained at point $\bx = (0, 0)$. 
Moreover, it shows that the errors of $\ap = 1.9$ are smaller than those of $\ap = 0.7$. 
In Table \ref{Table3-2-2}, we further explore numerical errors $\|e_u\|_\infty$ and convergence rates of our method for different $\ap$. 
It shows that our method with linear basis function has the second order of accuracy for any $\ap \in (0, 2)$, but  the smaller the power $\ap$, the less the numerical errors. 
 \begin{table}[htb!]
        \centering
        \begin{tabular}{|c||c|c|c|c|c|c|}
        \hline
        $h$ & 1/4 & 1/8 & 1/16 & 1/32 & 1/64 & 1/128 \\
        \hline
        \multirow{2}{*}{$\ap = 0.2$}&5.3181e-4 &  1.1946e-4  & 2.8883e-5 &  7.1505e-6  & 1.7827e-6 &  4.4531e-7\\
        \cline{2-7}
        & \mbox{c.r.}&2.1544 &  2.0482  & 2.0141  & 2.0040  & 2.0011\\
        \hline
        \multirow{2}{*}{$\ap = 0.7$}&2.3855e-3 &  5.1805e-4  & 1.2151e-4 &  2.9565e-5  & 7.3092e-6  & 1.8190e-6\\
        \cline{2-7}
        & \mbox{c.r.}&2.2031 &  2.0921 &  2.0391 &  2.0161 &  2.0065\\
        \hline
        \multirow{2}{*}{$\ap = 1$}&3.9406e-3  & 8.3747e-4 &  1.9083e-4 &  4.5384e-5  & 1.1056e-5  & 2.7276e-6\\
        \cline{2-7}
        & \mbox{c.r.}&2.2343 &  2.1338  & 2.0720 &  2.0374  & 2.0191\\
        \hline
        \multirow{2}{*}{$\ap = 1.4$}&6.9983e-3  & 1.4880e-3 &  3.2910e-4  & 7.5175e-5 &  1.7618e-5  & 4.2102e-6\\
        \cline{2-7}
        & \mbox{c.r.}&2.2337 &  2.1767  & 2.1302  & 2.0932 &  2.0651\\
        \hline
        \multirow{2}{*}{$\ap = 1.9$}&1.4264e-2  & 3.3824e-3  & 8.0943e-4 &  1.9424e-4  & 4.6676e-5  & 1.1230e-5\\
        \cline{2-7}
        & \mbox{c.r.}&2.0762 &  2.0631  & 2.0591 &  2.0571  & 2.0554\\
        \hline
        \end{tabular}
        \caption{Numerical errors $\|e_u\|_\infty$ and convergence rate (c.r.) in solving the 2D Poisson problem on $\Og = (-1, 1)^2$ with $f$ and $g$ defined in (\ref{fun1}), where linear basis $\varphi^1$ is used.}\label{Table3-2-2}
    \end{table}    
Extensive studies show that the performance of our method in two-dimensional cases are similar to those in one dimension, and we will omit showing details for brevity. 

\section{Conclusions and discussion}
\label{section4}
\setcounter{equation}{0}

We proposed a new class of operator factorization methods for discretizing the integral  fractional Laplacian $(-\Dt)^\fl{\ap}{2}$ in \eqref{fL-nD}. 
The performance of our method in both approximating the fractional Laplacian and solving the fractional Poisson problems was detailedly examined. 
The differentiation matrix  resulting from our method is of symmetric (multilevel) Toeplitz structure, which could save the memory cost and also enable efficient implementation via the fast Fourier transforms. 
Compared to other existing methods, our method has flexibility to increase numerical accuracy by using high-degree basis functions. 
Moreover,  the application of different basis functions only changes the entries of differentiation matrix, but does not affect its (multilevel) Toeplitz structure.  

Numerical studies show that for sufficiently smooth functions, our method with constant basis $\varphi^0$ and linear basis $\varphi^1$ has an optimal accuracy of ${\mathcal O}(h^2)$, and this rate can be improved to ${\mathcal O}(h^4)$ if quadratic basis $\varphi^2$ is used. 
The numerical accuracy can be further increased if higher degree basis function $\varphi^p$ is used. 
When approximating operator $(-\Dt)^\fl{\ap}{2}$,  the minimum condition to achieve the above optimal accuracy is  $u \in C^{2+\lfloor\ap\rfloor,\,\ap-\lfloor\ap\rfloor+\veps}(\bar{\Og})$ for constant and linear basis functions, while $u \in C^{4+\lfloor\ap\rfloor,\,\ap-\lfloor\ap\rfloor+\veps}(\bar{\Og})$ for quadratic basis for small $\veps > 0$.  
While solving fractional Poisson problems with  solution $u \in C^{m, l}(\bar{\Og})$ for $m \in {\mathbb N}$ and $0 < l < 1$, our method has an accuracy of ${\mathcal O}\big(h^{\min\{m+l,\,2\}}\big)$ for constant and linear basis functions, while ${\mathcal O}\big(h^{\min\{m+l,\,4\}}\big)$ for quadratic basis functions. 
Numerical experiments showed that choosing splitting parameter $\gm = 2$ leads to the best performance of our method, and thus this optimal parameter $\gm = 2$ should be always taken in our method. 
Our additional study on the tempered fractional Poisson problem confirmed that our method can be readily applied to a broader class of nonlocal operators.   

\bigskip
\noindent{\bf Acknowledgements. }  The authors thank Dr. Siwei Duo for her discussion and contributions at the initial stage of this project. 
This work was partially supported by the US National Science Foundation under Grant Number DMS-1913293 and DMS-1953177. 


\appendix
\section{Appendix}
\label{appendix1}

In one-dimensional cases, the weight integral in \eqref{wi-1D} can be found analytically. 
For convenience of the readers, we will summarize their analytical results of $\og_k^p$ for $p = 0, 1, 2$ as below.\\

For $p = 0$, the weight integral $\og_k^0$ in \eqref{wi-1D} is calculated as
\beas
\og_k^0 = \fl{1}{\sigma_0} \left\{\begin{array}{ll}
\xi_{1/2}^{\sigma_0},  & \mbox{if} \ \ k = 0,\\
\xi_N^{\sigma_0} - \xi_{N - 1/2}^{\sigma_0}, & \mbox{if} \ \ k = N,\\
\xi_{k + 1/2}^{\sigma_0} - \xi_{k - 1/2}^{\sigma_0}, \qquad\quad  & \mbox{otherwise},
\end{array}\right.
\eeas
where we denote $\sigma_0 = \gm - \ap$.  
For $p = 1$, the weight integral $\og_k^1$ in \eqref{wi-1D} is calculated as
\beas
\og_k^1 = \fl{h^{-1}}{\sigma_1\sigma_0} \left\{\begin{array}{ll}
\displaystyle h^{\sigma_1},  & \mbox{if} \ \ k = 0,\\
\displaystyle\xi_{N-1}^{\sigma_1} - \xi_{N}^{\sigma_1} + \sigma_1 h\,\xi_{N}^{\sigma_0}, & \mbox{if} \ \ k = N,\\
\displaystyle\xi_{k + 1}^{\sigma_1} + \xi_{k - 1}^{\sigma_1} - 2\xi_k^{\sigma_1}, \qquad\quad  & \mbox{otherwise},
\end{array}\right.
\eeas
where we denote $\sigma_1 = \sigma_0 + 1$. While  $p = 2$, the weight integral $\og_k^2$ in \eqref{wi-1D} is calculated as
\beas
\og_k^2 = \fl{1}{2h^2} \left\{\begin{array}{ll}
\displaystyle \fl{\xi_2^{\sigma_2}}{\sigma_2} - \xi_3 \fl{\xi_2^{\sigma_1}}{\sigma_1} + 2h^2\fl{\xi_2^{\sigma_0}}{\sigma_0},  & \mbox{if} \ \ k = 0,\\
&\\
\displaystyle \fl{\xi_N^{\sigma_2} - \xi_{N-2}^{\sigma_2}}{\sigma_2} - 2\xi_{N-\fl{3}{4}}\fl{\xi_N^{\sigma_1} - \xi_{N-2}^{\sigma_1}}{\sigma_1} + \xi_{N-1}\xi_{N-2}\fl{\xi_N^{\sigma_0} - \xi_{N-2}^{\sigma_0}}{\sigma_0}, & \mbox{if} \ \ k = N,\\
&\\
\displaystyle -2\fl{\xi_{k+1}^{\sigma_2} - \xi_{k-1}^{\sigma_2}}{\sigma_2} +4\xi_k\fl{\xi_{k+1}^{\sigma_1} - \xi_{k-1}^{\sigma_1}}{\sigma_1} -2\xi_{k+1}\xi_{k-1}\fl{\xi_{k+1}^{\sigma_0} - \xi_{k-1}^{\sigma_0}}{\sigma_0}, \qquad  & \mbox{if  \ $0 < k < N$,  \ odd}\\
&\\
\displaystyle \fl{\xi_{k+2}^{\sigma_2} - \xi_{k-2}^{\sigma_2}}{\sigma_2} - \fl{2}{\sigma_1}\big[\xi_{k-\fl{1}{2}}\big(\xi_{k}^{\sigma_1} - \xi_{k-2}^{\sigma_1}\big) +\xi_{k+\fl{1}{2}}\big(\xi_{k+2}^{\sigma_1} - \xi_{k}^{\sigma_1}\big)\big] & \\
\displaystyle \hspace{1cm}+ \fl{1}{\sigma_0}\big[\xi_{k-2}\xi_{k-1}\big(\xi_{k}^{\sigma_0} - \xi_{k-2}^{\sigma_0}\big)+\xi_{k+2}\xi_{k+1}\big(\xi_{k+2}^{\sigma_0} - \xi_{k}^{\sigma_0}\big)\big], & \mbox{if  \ $0 < k < N$,  \ even,}\\
\end{array}\right.
\eeas
where $\sigma_2 = \sigma_0 + 2$.

In high dimensions, the weight integral can not be calculated analytically, and thus numerical quadrature rules will be used to compute these integrals.  



\end{document}